\documentclass[leqno,11pt]{amsart}
\usepackage{amsmath,amscd}
\usepackage{epsfig,graphics}
\usepackage{amssymb,latexsym}
\usepackage{mathrsfs}
\newtheorem{thm}{\bf Theorem}[section]
\newtheorem{df}[thm]{\bf Definition}
\newtheorem{prop}[thm]{\bf Proposition}
\newtheorem{cor}[thm]{\bf Corollary}
\newtheorem{lem}[thm]{\bf Lemma}
\newtheorem{rem}[thm]{\bf Remark}
\newtheorem{ex}[thm]{\bf Example}

\newcommand{\A}{\mathcal{A}}
\newcommand{\B}{\mathcal{B}}

\newcommand{\cP}{\mathscr{P}}
\newcommand{\R}{\mathscr{R}}
\newcommand{\pf}{\noindent{\bfseries Proof. }}

\numberwithin{equation}{section}

\begin{document}
\title[Rational semistandard tableaux]
{Rational semistandard tableaux and character formula for the Lie
superalgebra  $\widehat{\frak{gl}}_{\infty|\infty}$}
\author{JAE-HOON KWON}
\address{Department of Mathematics \\ University of Seoul \\ 90
Cheonnong-dong, Dongdaemun-gu \\ Seoul 130-743, Korea }
\email{jhkwon@uos.ac.kr }

\thanks{This research was supported by KRF Grant $\sharp$2005-070-C00004 } 

\maketitle

\begin{abstract}
A new combinatorial interpretation of the Howe dual pair
$(\widehat{\frak{gl}}_{\infty|\infty},\frak{gl}_n)$ acting on an
infinite dimensional Fock space $\frak{F}^n$ of level $n$ is
presented. The character of a quasi-finite irreducible highest
weight representation of $\widehat{\frak{gl}}_{\infty|\infty}$
occurring in $\frak{F}^n$ is realized in terms of certain bitableaux
of skew shapes. We study a general combinatorics of these
bitableaux, including Robinson-Schensted-Knuth correspondence and
Littlewood-Richardson rule, and then its dual relation with the
rational semistandard tableaux for $\frak{gl}_n$. This result also
explains other Howe dual pairs including $\frak{gl}_n$.
\end{abstract}

\section{Introduction}
The Lie superalgebras and their representations appear naturally as
the fundamental algebraic structures in various areas of mathematics
and mathematical physics, and they have been studied by many people
since the fundamental work of Kac \cite{Kac77}.

Recently, from the viewpoint of Howe duality \cite{H,H92}, a close
relation between the representations of Lie algebras and Lie
superalgebras has been observed (see \cite{CW01,CW03,LS,N} and other
references therein), by which various character formulas for certain
Lie superalgebras have been obtained \cite{CL,CLZ,CZ}.

The purpose of this paper is to give a unified combinatorial
interpretation of the Howe dualities of the pairs
$(\frak{g},\frak{gl}_n)$ including a general linear Lie algebra.
Especially, our work will be devoted to the cases when
$\frak{g}=\widehat{\frak{gl}}_{\infty|\infty}$ or
$\widehat{\frak{gl}}_{\infty}$ acting on an infinite dimensional
Fock space which have been studied in
\cite{CL,CW01,CW03,KacR1,KacR2}. Let $\frak{F}^n$ ($n\geq 1$) be the
infinite dimensional Fock space generated by $n$ pairs of free
fermions and $n$ pairs of free bosons (see \cite{CL,CW03}). On
$\frak{F}^n$, there exists a natural commuting actions of the
infinite dimensional Lie superalgebra
$\widehat{\frak{gl}}_{\infty|\infty}$ and the finite dimensional Lie
algebra $\frak{gl}_n$. Using Howe duality, Cheng and Wang derived a
multiplicity-free decomposition
\begin{equation*}
\frak{F}^n\simeq \bigoplus_{\lambda\in
\mathbb{Z}_+^n}L_{\lambda}\otimes L_n(\lambda),
\end{equation*}
as a $\widehat{\frak{gl}}_{\infty|\infty}\oplus \frak{gl}_n$-module,
where the sum ranges over all generalized partitions $\lambda$ of
length $n$, $L_n(\lambda)$ is the rational representation of
$\frak{gl}_n$ corresponding to $\lambda$, and $L_{\lambda}$ is the
associated quasi-finite irreducible highest weight representation of
$\widehat{\frak{gl}}_{\infty|\infty}$ \cite{CW03}. From the above
decomposition and the classical Cauchy identities for hook Schur
polynomials (cf.\cite{BR,Rem}), Cheng and Lam derived a character
formula of $L_{\lambda}$, in terms of hook Schur polynomials
\cite{CL}. They also described the tensor product decomposition of
$L_{\lambda}\otimes L_{\mu}$ for $\lambda\in\mathbb{Z}_+^m$ and
$\mu\in\mathbb{Z}_+^n$. In fact, this is a natural super-analogue of
the Kac and Radul's works on the construction of quasi-finite
irreducible highest weight representation $L^0_{\lambda}$ of the
infinite dimensional Lie algebra $\widehat{\frak{gl}}_{\infty}$,
which is also parameterized by a generalized partition \cite{KacR2}.

Suppose that $\A$ and $\B$ are linearly ordered
$\mathbb{Z}_2$-graded sets. Motivated by the character formula  of
$L_{\lambda}$ or $L^0_{\lambda}$ in \cite{CL,KacR2}, we introduce
the notion of {\it $\A/\B$-semistandard tableaux of shape
$\lambda$}, which is our main object in this paper. Roughly
speaking, an $\A/\B$-semistandard tableau of shape $\lambda$ is a
pair of tableaux $(T^+,T^-)$ such that $T^+$ (resp. $T^-$) is a
semistandard tableau with letters in $\A$ (resp. $\B$), where the
shapes of $T^+$ and $T^-$ are not necessarily fixed ones but satisfy
certain conditions determined by $\lambda$. We develop the insertion
scheme for $\A/\B$-semistandard tableaux, and derive analogues of
Robinson-Schensted-Knuth (or simply RSK) correspondence and
Littlewood-Richardson (or simply LR) rule.

As in the case of ordinary semistandard tableaux, we define {\it
skew $\A/\B$-semistandard tableaux of shape $\lambda/\mu$} for
arbitrary two generalized partitions $\lambda$ and $\mu$ of the same
length, and describe the corresponding skew LR rule. Then it turns
out that the combinatorics of $\A/\B$-semistandard tableaux is {\it
dual} to that of rational semistandard tableaux for general linear
Lie algebra introduced by Stembridge \cite{St} in the sense that the
skew LR rule (resp. LR rule) of $\A/\B$-semistandard tableaux are
completely determined by the LR rule (resp. skew LR rule) of
rational semistandard tableaux.

Next, we show that the character of $SST_{\A/\B}(\lambda)$, the set
of all $\A/\B$-semistandard tableaux of shape $\lambda$, reduce to
the character of $L_{\lambda}$ or $L^0_{\lambda}$ under suitable
choices of $\A$ and $\B$. This is done by observing that we have
another expression of the character of $SST_{\A/\B}(\lambda)$, which
is a kind of branching formula very similar to the Cheng and Lam's
formula (or the Kac and Radul's formula). As a corollary, we
immediately obtain new combinatorial interpretations of the
decomposition of the Fock space representation $\frak{F}^n$ and the
tensor product $L_{\lambda}\otimes L_{\mu}$ from RSK correspondence
and LR rule, respectively. We also obtain a Jacobi-Trudi type
character formula for $L_{\lambda}$ or $L^0_{\lambda}$, which has
not been observed yet as far as we know. The dual relationship
between $\A/\B$-semistandard tableaux and rational semistandard
tableaux now explains the duality between  the tensor product
decomposition and the branching rule of the pairs
$(\widehat{\frak{gl}}_{\infty|\infty},\frak{gl}_n)$ and
$(\widehat{\frak{gl}}_{\infty},\frak{gl}_n)$, which is a general
feature in any Howe dual pair. We expect a combinatorial
construction of $L_{\lambda}$ as a vector space spanned by
$SST_{\A/\B}(\lambda)$, and also a $q$-analogue  of $L_{\lambda}$ as
a representation of the associated quantum group
$U_q(\widehat{\frak{gl}}_{\infty|\infty})$ with a crystal graph
$SST_{\A/\B}(\lambda)$.

Finally, we would like to remark that the notion of
$\A/\B$-semistandard tableaux can be applied to other classes of
representations of Lie (super)algebras. For example, when $\A$ is
finite and $\B$ is empty, our results explain the classical Howe
duality of the $(\frak{gl}_{p|q},\frak{gl}_n)$ pair acting on
$S(\mathbb{C}^{p|q}\otimes \mathbb{C}^n)$, the supersymmetric
algebra generated by $\mathbb{C}^{p|q}\otimes \mathbb{C}^n$
(cf.\cite{CW01,K}). Moreover, when $\A$ and $\B$ are both non-empty
finite sets, we can recover the combinatorial picture of the Howe
dual pair $(\frak{gl}_{p|q},\frak{gl}_n)$ on a supersymmetric
algebra \cite{CLZ} (see also \cite{KV}), and $SST_{\A/\B}(\lambda)$
($\lambda\in\mathbb{Z}_+^n$) realizes the character of an
infinite-dimensional unitarizable representation of
$\frak{gl}_{p|q}$. In general, we expect that to arbitrary $\A$ and
$\B$, there corresponds a contragredient Lie superalgebra, where the
character of $SST_{\A/\B}(\lambda)$ gives the character of an
irreducible representation parameterized by $\lambda$ (cf.\cite{K}).
\vskip 3mm

%We can say that the Grothendieck ring for the category of
%semi-simple representations of $\widehat{\frak{gl}}_{\infty|\infty}$
%(or $\widehat{\frak{gl}}_{\infty}$) whose irreducible factors ranges
%over $L_{\lambda}$'s (or $L^0_{\lambda}$'s) with finite multiplicity
%for each $\lambda$, is isomorphic as a commutative
%$\mathbb{Z}$-algebra to the dual of the coalgebra of the graded sum
%of Grothendieck groups for the category of rational representations
%of $\frak{gl}_n$, whose comultiplication is induced from the
%canonical branching rule.

This paper is organized as follows. In Section 2, we briefly recall
a necessary background on semistandard tableaux for Lie
superalgebras and their insertion scheme  including LR rule. In
Section 3, we review the notion of rational semistandard tableaux
and define $\A/\B$-semistandard tableaux. Some properties of the
characters are also discussed. In Section 4, we introduce an
insertion algorithm for $\A/\B$-semistandard tableaux, and then
derive analogues of RSK correspondence, and (skew) LR rule. Finally,
in Section 5, we show that the character of $SST_{\A/\B}(\lambda)$
reduces to the character of $L_{\lambda}$ or $L^0_{\lambda}$ under a
particular choice of $\A$ and $\B$, and then explain the
relationship between the combinatorial results established in the
previous sections and the representations of
$\widehat{\frak{gl}}_{\infty|\infty}$ or
$\widehat{\frak{gl}}_{\infty}$.\vskip 3mm

{\bf Acknowledgement} This work was done during the author's visit
at Korea Institute for Advanced Study in 2006. He thanks the
institute for its hospitality. He is also grateful to Prof. S.-J.
Cheng for helpful comments.

\section{Preliminaries}
Throughout the paper, we assume that $\A$ and $\B$ are linearly
ordered sets which are at most countable, and also
$\mathbb{Z}_2$-graded (that is, $\A=\A_0\sqcup\A_1$ and
$\B=\B_0\sqcup\B_1$). By convention, we let
$\mathbb{N}=\{\,1<2<\cdots\,\}$, and $[n]=\{\,1<\cdots<n\,\}$ for
$n\geq 1$, where all the elements are of degree $0$. \vskip 3mm

In this section, we introduce the notion of $\A$-semistandard
tableaux, and describe the associated Littlewood-Richardson rule. In
fact, the combinatorics of $\A$-semistandard tableaux is essentially
the same as that of semistandard tableaux with entries in
$\mathbb{N}$. Since the only difference is that we assume the {\it
column strict} condition on entries of degree $0$ and the {\it row
strict} condition on entries of degree $1$, most of the results in
this section can be verified directly by modifying the arguments of
the corresponding results in case of ordinary semistandard tableaux
(cf.\cite{Fu,KK,Rem}). So, we leave the detailed verifications to
the readers.

\subsection{Semistandard tableaux}
Let us recall some basic terminologies (cf.\cite{Mac}). A {\it
partition} of a non-negative integer $n$ is a non-increasing
sequence of non-negative integers $\lambda = (\lambda_k)_{k\geq 1}$
such that $\sum_{k\geq 1}\lambda_k=n$. We also write $|\lambda|=n$.
Each $\lambda_k$ is called a {\it part of $\lambda$}, and the number
of non-zero parts is called the {\it length of $\lambda$} denoted as
$\ell(\lambda)$. We also write $\lambda=(1^{m_1},2^{m_2},\cdots)$,
where $m_i$ is the number of occurrences of $i$ in $\lambda$. We
denote by $\cP$ the set of all partitions, and denote by
$\cP_{\ell}$ ($\ell\geq 1$) the set of all partitions with length no
more than $\ell$. Recall that a partition $\lambda =
(\lambda_k)_{k\geq 1}$ is identified with a {\it Young diagram}
which is a collection of nodes (or boxes) in left-justified rows
with $\lambda_k$ nodes in the $k$th row. For $\lambda\in\cP$,
$\lambda'$ is the conjugate of $\lambda$. For $\lambda,\mu\in\cP$
with $\lambda\supset\mu$ (that is, $\lambda_i\geq \mu_i$ for all
$i$), $\lambda/\mu$ is the skew Young diagram obtained from
$\lambda$ by removing $\mu$, and $|\lambda/\mu|$ is defined to be
the number of nodes in $\lambda/\mu$.

\begin{df}\label{A SST}{\rm
For a skew Young diagram $\lambda/\mu$, a tableau $T$ obtained by
filling $\lambda/\mu$ with entries in $\A$ is called {\it
$\A$-semistandard} if
\begin{itemize}
\item[(1)] the entries in each row (resp. column) are weakly increasing from left to right (resp. from top to bottom),
\item[(2)] the entries in $\A_0$ (resp. $\A_1$) are strictly increasing in each column (resp.
row).
\end{itemize}  }
\end{df}
We say that $\lambda/\mu$ is the {\it shape of $T$}, and write ${\rm
sh}(T)=\lambda/\mu$. We denote by $SST_{\A}(\lambda/\mu)$ the set of
all $\A$-semistandard tableaux of shape $\lambda/\mu$. We assume
that $SST_{\A}(\lambda)=\{\emptyset\}$ when $\lambda\in\cP$ is a
partition of $0$, where $\emptyset$ denotes the empty tableau. We
set $\cP_{\A}=\{\,\lambda\in\cP\,|\,SST_{\A}(\lambda)\neq
\emptyset\,\}$. Note that $\cP_{[n]}=\cP_n$, and $\cP_{\A}=\cP$ when
$\A$ is an infinite set.

For $T\in SST_{\A}(\lambda/\mu)$, we let
$$w_{\rm col}(T)=a_1a_2\ldots a_n$$ be the word with letters in
$\A$ obtained by reading the entries of $T$ column by column from
right to left, and from top to bottom in each column. Also, we let
$$w_{\rm row}(T)=b_1b_2\ldots b_n$$ be the word
obtained by reading the entries of $T$ row by row from bottom to
top, and from left to right in each row.\vskip 5mm

\subsection{Operations on tableaux} Let us define several operations on semistandard
tableaux. Let $\lambda/\mu$ be a skew Young diagram.

\begin{itemize}
\item[(1)] {\it transpose} :   We define $\A'$ to be the linearly ordered
$\mathbb{Z}_2$-graded set such that $\A'=\A$ with the same linear
ordering as in $\A$ and the opposite $\mathbb{Z}_2$-grading, that
is, $\A'_0=\A_1$, $\A'_1=\A_0$. For $T\in SST_{\A}(\lambda/\mu)$, we
denote by $T^t$ the transpose of $T$. Then  $T^t\in
SST_{\A'}(\lambda'/\mu')$.

\item[(2)] {\it rotation} :  We define $\A^{\pi}$ to be the linearly ordered
$\mathbb{Z}_2$-graded set such that $\A^{\pi}=\A$ with the same
$\mathbb{Z}_2$-grading and the reverse linear ordering. For $T\in
SST_{\A}(\lambda/\mu)$, we define $T^{\pi}$ to be the tableau
obtained by applying $180^{\circ}$-rotation to $T$. Then
$T^{\pi}\in SST_{\A^{\pi}}((\lambda/\mu)^{\pi})$, where
$(\lambda/\mu)^{\pi}$ is the skew diagram obtained from
$\lambda/\mu$ after $180^{\circ}$-rotation.

\item[(3)] {\it reverse transpose} : We set $\A^{\sharp}=(\A')^{\pi}$. For $T\in
SST_{\A}(\lambda/\mu)$, we define $T^{\sharp}=(T^t)^{\pi}$. Then
$T^{\sharp}\in SST_{\A^{\sharp}}((\lambda/\mu)^{\sharp})$, where
$(\lambda/\mu)^{\sharp}=(\lambda'/\mu')^{\pi}$.
\end{itemize}

We define $\A\ast\B$ to be the $\mathbb{Z}_2$-graded set $\A\sqcup
\B$ with the extended linearly ordering given by $a<b$ for all
$a\in\A$ and $b\in\B$. For $S\in SST_{\A}(\mu)$ and $T\in
SST_{\B}(\lambda/\mu)$, we define $S\ast T$ to be the
$\A\ast\B$-semistandard tableau of shape $\lambda$ obtained by
gluing $S$ and $T$ so that the right-most node in each row of $S$ is
placed next to the left-most node in the same row of $T$.\vskip 5mm

\subsection{Insertions} Let us describe the {\it Schensted's column and row bumping
algorithms} for $\A$-semistandard tableaux (cf.\cite{BR,Rem}):
Suppose that $a\in\A$ and $T\in SST_{\A}(\lambda)$
($\lambda\in\cP_{\A}$) are given.\vskip 3mm

First, we define $(T\leftarrow a )$ to be the tableau obtained
from $T$ by applying the following procedure (called the {\it
column bumping algorithm});
\begin{itemize}
\item[(i)] If $a\in \A_0$, let $a'$ be the smallest entry in
the first (or the left-most) column, which is greater than or
equal to $a$. If $a\in \A_1$, let $a'$ be the smallest entry in
the first column, which is greater than $a$. If there are more
than one $a'$, choose the one in the highest position.

\item[(ii)]  Replace {$a'$} by {$a$}. If there is no such $a'$, put
{$a$} at the bottom of the first column and stop the procedure.

\item[(iii)]  Repeat (i) and (ii) on the next column with {$a'$}.
\end{itemize}\vskip 3mm

Next, we define $(a \rightarrow T )$  to be the tableau obtained
from $T$ by applying the following procedure (called the {\it row
bumping algorithm});
\begin{itemize}
\item[(i)]  If $a\in \A_0$, let $a'$ be the smallest entry in the
first (or the top) row, which is greater than $a$.  If $a\in
\A_1$, let $a'$ be the smallest entry in the first row, which is
greater than or equal to $a$. If there are more than one $a'$,
choose the one in the left-most position.

\item[(ii)] Replace {$a'$} by {$a$}. If there is no such $a'$, put
{$a$} at the right-most end of the top row and stop the procedure.

\item[(iii)]  Repeat (i) and (ii) on the next row with {$a'$}.
\end{itemize}\vskip 3mm

Suppose that $\mu,\nu\in \cP_{\A}$ are given. For $T\in
SST_{\A}(\mu)$ and $T'\in SST_{\A}(\nu)$, let $w_{\rm
col}(T')=c_1c_2\ldots c_n$. We define
\begin{equation*}
(T\leftarrow T')=((((T\leftarrow c_1)\leftarrow c_2)\cdots
)\leftarrow c_n).
\end{equation*}
Similarly, let $w_{\rm row}(T')=r_1r_2\ldots r_m$. We define
\begin{equation*}
(T'\rightarrow T)=(
r_m\rightarrow(\cdots(r_2\rightarrow(r_1\rightarrow T)))).
\end{equation*}

We define $(T\leftarrow T')_R$ to be the semistandard tableau in
$SST_{\mathbb{N}}(\lambda/\mu)$ ($\lambda={\rm sh}(T\leftarrow T')$)
such that if ${c_i}$ is in the $k$th row of $T'$ and inserted into
$(((\cdots(\,T\leftarrow {c_1}\,)\leftarrow c_2)\cdots)\leftarrow
{c_{i-1}})$ to create a node in $\lambda/\mu$, then we fill the node
with $k$. We call $(T\leftarrow T')_R$ the {\it recording tableau of
$(T\leftarrow T')$}.

Similarly, we define $(T'\rightarrow T)_R$ to be the semistandard
tableau in $SST_{\mathbb{N}'}(\eta/\mu)$ ($\eta={\rm
sh}(T'\rightarrow T)$) such that if ${r_i}$ is in the $k$th column
of $T'$ and inserted into
$(r_{i-1}\rightarrow(\cdots(r_2\rightarrow(r_1\rightarrow T))))$ to
create a node in $\eta/\mu$, then we fill the node with  $k$. We
also call $(T'\rightarrow T)_R$ the {\it recording tableau of
$(T'\rightarrow T)$}.

Given $\lambda$, $\mu$ and $\nu$ in $\cP$ such that $\mu\subset
\lambda$ and $|\lambda|=|\mu|+|\nu|$, a tableau $T\in
SST_{\mathbb{N}}(\lambda/\mu)$ is called a {\it
Littlewood-Richardson tableau  of shape $\lambda/\mu$ with content
$\nu$} if
\begin{itemize}
\item[(1)] the number of occurrences of $k$ in $T$ is equal to
$\mu_k$ for $k\geq 1$,

\item[(2)] $w_{\rm col}(T)$ is a lattice permutation (see
\cite{Mac} for its definition).
\end{itemize}
We denote by $LR^{\lambda}_{\mu\,\nu}$ the set of all
Littlewood-Richardson tableaux of shape $\lambda/\mu$ with content
$\nu$, and put $|LR^{\lambda}_{\mu\,\nu}|=N^{\lambda}_{\mu\,
\nu}$, which is called  a {\it Littlewood-Richardson coefficient}.
By standard arguments as in the case of $\mathbb{N}$-semistandard
tableaux, we can check the following.

\begin{lem}[cf.\cite{KK,Rem,Thomas}]\label{LR tableaux}
For $T\in SST_{\A}(\mu),\ T'\in SST_{\A}(\nu)$
$(\mu,\nu\in\cP_{\A})$,
\begin{itemize}
\item[(1)] $(T\leftarrow T')_R \in LR^{\lambda}_{\mu\,\nu}$, where
$\lambda={\rm sh}(T\leftarrow T')$,
\item[(2)] $[(T'\rightarrow T)_R]^t \in LR^{\eta'}_{\mu'\,\nu'}$,
where $\eta={\rm sh}(T'\rightarrow T)$.
\end{itemize} \qed
\end{lem}\vskip 3mm

By Lemma \ref{LR tableaux}, we obtain the {\it Littlewood-Richardson
rule} (or simply LR rule) for $\A$-semistandard tableaux.
\begin{thm}[cf.\cite{KK,Rem,Thomas}]\label{LR rule} Suppose that $\mu,\nu\in\cP_{\A}$ are given.
\begin{itemize}
\item[(1)] The map $\rho_{\rm col} : (T,T')\mapsto ((T\leftarrow T'),(T\leftarrow
T')_R)$ gives a bijection
$$\rho_{\rm col} : SST_{\A}(\mu)\times SST_{\A}(\nu)\longrightarrow
\bigsqcup_{\lambda\in\cP_{\A}}SST_{\A}(\lambda)\times
LR^{\lambda}_{\mu\,\nu}.$$

\item[(2)] The map $\rho_{\rm row} : (T,T')\mapsto ((T'\rightarrow T),[(T'\rightarrow T)_R]^t)$ gives a bijection
$$\rho_{\rm row} : SST_{\A}(\mu)\times SST_{\A}(\nu)\longrightarrow
\bigsqcup_{\eta\in\cP_{\A}}SST_{\A}(\eta)\times LR^{\eta'}_{\mu'\,
\nu'}.$$
\end{itemize}\qed
\end{thm}

For an $r$-tuple of non-negative integers $\nu=(\nu_1,\cdots,\nu_r)$
such that $\nu_i\in\cP_{\A}$ ($1\leq i\leq r$), consider
\begin{equation*}
(T_1,\ldots,T_r)\in SST_{\A}(\nu_1)\times\cdots\times
SST_{\A}(\nu_r).
\end{equation*}
For $1\leq i\leq r$, put
\begin{equation*}
S_i=((((T_1\leftarrow T_2)\leftarrow T_3)\cdots )\leftarrow T_i),
\end{equation*}
and $S=S_r$. If we put $\mu^{(i)}={\rm sh}(S_i)\in \cP$, then we
have $\nu_1=\mu^{(1)}\subset\cdots\subset \mu^{(r)}=\mu$, and
$\mu^{(i)}/\mu^{(i-1)}$ ($1\leq i\leq r$) is a horizontal strip of
length $\nu_i$ (we assume that $\mu^{(0)}$ is the empty partition).
Filling  $\mu^{(i)}/\mu^{(i-1)}$ with $i$, we obtain an
$[r]$-semistandard tableau $S_R\in SST_{[r]}(\mu)$, with content
$\nu$ (that is, each entry $i$ occurs as many times as $\nu_i$ for
$1\leq i\leq r$). The correspondence $(T_1,\ldots,T_r)\mapsto
(S,S_R)$ is reversible by Theorem \ref{LR rule}.

Similarly, put
\begin{equation*}
S'_i=( T_i\rightarrow(\cdots(T_3\rightarrow(T_2\rightarrow T_1)))),
\end{equation*}
for $1\leq i\leq r$ and $S'=S'_r$  . If we put $\mu^{(i)}={\rm
sh}(S'_i)\in \cP$, then we have $\nu_1=\mu^{(1)}\subset\cdots\subset
\mu^{(r)}=\mu$, and $\mu^{(i)}/\mu^{(i-1)}$ ($1\leq i\leq r$) is
also a horizontal strip of length $\nu_i$. So, as in the case of
$S_R$, we may define an $[r]$-semistandard tableau of shape $\mu$
with content $\nu$, say $S'_R$. The correspondence
$(T_1,\ldots,T_r)\mapsto (S',S'_R)$ is also reversible.

Summarizing the arguments, we have
\begin{prop}\label{row insertions} Under the above hypothesis, we
have two bijections
\begin{equation*}
\varrho_{\rm col}, \varrho_{\rm row} :
SST_{\A}(\nu_1)\times\cdots\times SST_{\A}(\nu_r) \longrightarrow
\bigsqcup_{\mu\in\cP_{\A}}SST_{\A}(\mu)\times SST_{[r]}(\mu)_{\nu},
\end{equation*}
where $\varrho_{\rm col}(T_1,\ldots,T_r)=(S,S_R)$, $\varrho_{\rm
row}(T_1,\ldots,T_r)=(S',S'_R)$, and $SST_{[r]}(\mu)_{\nu}$ is the
set of all $[r]$-semistandard tableaux of shape $\mu$ with content
$\nu$.\qed
\end{prop}\vskip 5mm

\subsection{Switching algorithm} Let us describe the skew LR rule
for $\A$-semistandard tableaux. To do this, we will use the {\it
switching algorithm} introduced by Benkart, Sottile, and Stroomer
\cite{BSS}.

Let $\lambda/\mu$ be a skew Young diagram. Let $U$ be a tableau of
shape $\lambda/\mu$ with entries in $\A\sqcup\B$, satisfying the
following conditions;
\begin{itemize}
\item[(S1)] if $u,u'\in \A$ (resp. $\B$) are entries of $U$ and $u$ is northwest of $u'$,
then $u\leq u'$,

\item[(S2)] in each column of $U$, entries in $\A_0$ or $\B_0$ increase
strictly,

\item[(S3)] in each row of $U$, entries in $\A_1$ or $\B_1$ increase
strictly,
\end{itemize}
where we say that {\it $u$ is northwest of $u'$} provided the row
and column indices of $u$ are no more than those of $u'$.

Suppose that $a\in \A$ and $b\in \B$ are two adjacent entries in
$U$ such that $a$ is placed above or to the left of $b$.
Interchanging $a$ and $b$ is called a {\it switching} if the
resulting tableau still satisfies the conditions (S1), (S2) and
(S3).

\begin{thm}[Theorem 2.2 and 3.1 in \cite{BSS}]\label{switching}
Let $\lambda/\mu$ be a skew Young diagram. For $S\in SST_{\A}(\mu)$
and $T\in SST_{\B}(\lambda/\mu)$, let $U$ be a tableau obtained from
$S\ast T$ by applying switching procedures as far as possible. Then
\begin{itemize}
\item[(1)] $U=T'\ast S'$,
where $T'\in SST_{\B}(\nu)$ and $S'\in SST_{\A}(\lambda/\nu)$ for
some $\nu$.

\item[(2)] $U$ is uniquely determined by $S$ and $T$.

\item[(3)] When $\A=\mathbb{N}$, $S'\in
LR^{\lambda}_{\nu\,\mu}$ if and only if $S=H^{\mu}$, where $H^{\mu}$
is the unique $\mathbb{N}$-semistandard tableau of shape $\mu$ with
content $\mu$.
\end{itemize}\qed
\end{thm}

\begin{ex}{\rm
Suppose that $\A=\mathbb{N}$ and
$\B=\mathbb{N}'=\{1'<2'<3'<\ldots\}$ (see 2.2). Consider
$$S\ast T=
\begin{array}{ccc}
  1 & 1 & 2 \\
  2 & 3 & 3' \\
  1' & 2' & 3'
\end{array}\in SST_{\mathbb{N}\ast\mathbb{N}'}(3^3),
$$
where $S\in SST_{\mathbb{N}}(3,2)$, and $T\in
SST_{\mathbb{N}'}((3^3)/(3,2))$. Then
\begin{equation*}
\begin{split}
&S\ast T = \begin{array}{ccc}
  1 & 1 & 2 \\
  2 & 3 & 3' \\
  1' & 2' & 3'
\end{array} \stackrel{ 1'\leftrightarrow 2}{\longrightarrow}
\begin{array}{ccc}
  1 & 1 & 2 \\
  1' & 3 & 3' \\
  2 & 2' & 3'
\end{array} \stackrel{1'\leftrightarrow 1}{\longrightarrow}
\begin{array}{ccc}
  1' & 1 & 2 \\
  1 & 3 & 3' \\
  2 & 2' & 3'
\end{array} \\
& \stackrel{2'\leftrightarrow 3}{\longrightarrow}
\begin{array}{ccc}
  1' & 1 & 2 \\
  1 & 2' & 3' \\
  2 & 3 & 3'
\end{array} \stackrel{2'\leftrightarrow 1}{\longrightarrow}
\begin{array}{ccc}
  1' & 2' & 2 \\
  1 & 1 & 3' \\
  2 & 3 & 3'
\end{array} \stackrel{3'\leftrightarrow 2}{\longrightarrow}
\begin{array}{ccc}
  1' & 2' & 3' \\
  1 & 1 & 2 \\
  2 & 3 & 3'
\end{array} \\
& \stackrel{3'\leftrightarrow 3}{\longrightarrow}
\begin{array}{ccc}
  1' & 2' & 3' \\
  1 & 1 & 2 \\
  2 & 3' & 3
\end{array} \stackrel{3'\leftrightarrow 2}{\longrightarrow}
\begin{array}{ccc}
  1' & 2' & 3' \\
  1 & 1 & 2 \\
  3' & 2 & 3
\end{array} \stackrel{3'\leftrightarrow 1}{\longrightarrow}
\begin{array}{ccc}
  1' & 2' & 3' \\
  3' & 1 & 2 \\
  1 & 2 & 3
\end{array} = T' \ast S'.
\end{split}
\end{equation*}
}
\end{ex}\vskip 3mm

\begin{rem}{\rm (1)
Theorem 2.2 in \cite{BSS} is shown when $\A=\B=\mathbb{N}$ having
only elements of degree 0. But we may naturally extend this result
to arbitrary $\A$ and $\B$. We leave the detailed verifications to
the readers (see also \cite{Rem}).

(2) The resulting tableau $T'$ in Theorem \ref{switching} is
independent of the choice of $S$, and the algorithm of producing
$T'$ from $T$ is known as the {\it Sch\"{u}tzenberger's jeu de
taquin slides}. }
\end{rem}

Suppose that  $T\in SST_{\A}(\lambda/\mu)$ is given for a skew Young
diagram $\lambda/\mu$. Consider
\begin{equation*}
H^{\mu}\ast T \in SST_{\mathbb{N}\ast\A}(\lambda).
\end{equation*}
Using the switching procedures in Theorem \ref{switching}, we obtain
a unique $\A\ast\mathbb{N}$-semistandard tableau given by
\begin{equation*}
T'\ast U \in SST_{\A\ast\mathbb{N}}(\lambda),
\end{equation*}
where $T'\in SST_{\A}(\nu)$ for some $\nu$, and $U\in
LR^{\lambda}_{\nu\,\mu}$ (see also Example 3.3 in \cite{BSS}). Let
us put
$$j(T)=T', \ \ \ \ \ j(T)_R=U.$$

\begin{cor}[\cite{BSS}]\label{skew LR} Under the above hypothesis,
\begin{itemize}
\item[(1)] the map $J : T\mapsto (j(T),j(T)_R)$ gives a
bijection $$J : SST_{\A}(\lambda/\mu) \longrightarrow
\bigsqcup_{\nu\in\cP_{\A}} SST_{\A}(\nu)\times
LR^{\lambda}_{\nu\,\mu},$$

\item[(2)] if $\A=\mathbb{N}$, then the map $Q \mapsto j(Q)_R$ restricts to a bijection from
$LR^{\lambda}_{\mu\,\nu}$ to $LR^{\lambda}_{\nu\,\mu}$. In
particular, we have $N^{\lambda}_{\mu\,\nu}=N^{\lambda}_{\nu\,\mu}$.
\end{itemize}\qed
\end{cor}
Let us remark another symmetry of Littlewood-Richardson tableaux,
a bijective proof of which has been given in \cite{HS}. For
self-containedness, we give another proof using switching
algorithm (see also \cite{BSS}).
\begin{cor}[cf.\cite{BSS,HS}]\label{conjugate LR} There exists a bijection
$\tau : LR^{\lambda}_{\mu\,\nu}\rightarrow
LR^{\lambda'}_{\mu'\,\nu'}$ for $\lambda,\mu,\nu\in\cP$.
\end{cor}
\pf Given $T\in LR^{\lambda}_{\mu\,\nu}$, consider $S\ast T^t$ where
$S=H^{\mu'}\in SST_{\mathbb{N}}(\mu')$. Note that $S\ast T^t$ is an
$\mathbb{N}\ast \mathbb{N}'$-semistandard tableau of shape
$\lambda'$. Applying the switching procedures to the pair $(S,T^t)$
as far as possible, we obtain by Theorem \ref{switching} (3) an
$\mathbb{N}'\ast \mathbb{N}$-semistandard tableau of shape
$\lambda'$
$$(H^{\nu})^t\ast S',$$
where $S'\in LR^{\lambda'}_{\nu'\,\mu'}$. Now, if we define
$\tau(T)=j(S')_R$ (see Corollary \ref{skew LR} (2)), then $\tau$ is
a bijection between $LR^{\lambda}_{\mu\,\nu}$ and
$LR^{\lambda'}_{\mu'\,\nu'}$. \qed\vskip 3mm

One may define semistandard tableaux with entries from a given
$\mathbb{Z}_2$-graded set $\A$ with different linear orderings $<$
and $<'$. Then the switching algorithm enables us to construct
easily a bijection between these two kinds of semistandard tableaux
of a given skew shape (cf.\cite{BSS,Hai,Rem}).

 First, consider $\A_0\ast\A_1$, where the linear
orderings on $\A_0$ and $\A_1$ remain the same. We may view
$\A_0\ast\A_1$ as a {\it shuffling of $\A$}.
\begin{prop}\label{shuffling}
For a skew Young diagram $\lambda/\mu$, there exists a bijection
between $SST_{\A}(\lambda/\mu)$ and
$SST_{\A_0\ast\A_1}(\lambda/\mu)$.
\end{prop}
\pf Let $<$ denotes the linear ordering on $\A$. Given $T\in
SST_{\A}(\lambda/\mu)$, suppose that $T$ is not
$\A_0\ast\A_1$-semistandard. Then there exists an entry $a$ of
degree 1 in $T$ such that $a<a'$ for some entry $a'$ of degree 0 in
$T$. Let $a_{\rm max}$ be the largest such entry of degree 1 in $T$.
Also, we let $a_{\rm min}$ be the smallest entry of degree 1 in $T$
such that there exists no entry of degree 0 in $T$ greater than
$a_{\rm min}$. Note that $a_{\rm max}<a_{\rm min}$. By definition,
we can check that there is no entry $a$ of degree 1 in $T$ such that
$a_{\rm max}<a<a_{\rm min}$.

Consider the tableau $S$ obtained by removing all the entries of $T$
smaller than $a_{\rm max}$ or no less than $a_{\rm min}$ (if there
is no such $a_{\rm min}$, then we assume $a_{\rm min}$ to be a
formal symbol greater than any element in $\A$). Then
\begin{equation*}
S=S'\ast S'',
\end{equation*}
where $S'$ is a tableaux of a skew shape with the entry $a_{\rm
max}$, and $S'$ is an $\A_0$-semistandard tableau with entries in
$\{\,a\in\A_0\,|\,a_{\rm max}<a<a_{\rm min}\,\}$. Now if we apply
the switching algorithm in Theorem \ref{switching} to the pair
$(S',S'')$, then we get a new tableau $U$ of the same shape as $S$
such that
\begin{equation*}
U=U'\ast U'',
\end{equation*}
where $U'$ is an $\A_0$-semistandard tableau with entries in
$\{\,a\in\A_0\,|\,a_{\rm max}<a<a_{\rm min}\,\}$, and $U''$ is a
tableaux of a skew shape with the entry $a_{\rm max}$.

Repeating the above argument, we obtain a unique tableau $T^*\in
SST_{\A_0\ast \A_1}(\lambda/\mu)$. By construction, this gives the
required bijective correspondence. \qed\vskip 3mm

By similar arguments, we also have
\begin{prop}\label{reverse ordering}
For a skew Young diagram $\lambda/\mu$, there exists a bijection
between $SST_{\A}(\lambda/\mu)$ and
$SST_{\A^{\pi}}(\lambda/\mu)$.\qed
\end{prop}\vskip 5mm

\section{$\A/\B$-semistandard tableaux}
\subsection{Rational semistandard tableaux}\label{rational tableaux}
First, let us recall the notion of rational $[n]$-semistandard
tableaux ($n\geq 1$) which was introduced by Stembridge \cite{St}.

For $n\geq 1$, let
\begin{equation}
\mathbb{Z}^n_+=\{\,\lambda=(\lambda_1,\cdots,\lambda_n)\in\mathbb{Z}^n\,|\,\lambda_1\geq
\lambda_2\geq\cdots\geq \lambda_n\,\}
\end{equation}
be the set of all {\it generalized partitions of level $n$} (we will
call $n$ as the level of $\lambda$ rather than its length to avoid a
confusion with the length of a partition). We have a natural
embedding $\cP_n\subset \mathbb{Z}_+^n$. For
$\lambda\in\mathbb{Z}_+^n$, we write $\langle\lambda
\rangle=\sum_{1\leq i\leq n}\lambda_i$, and $|\lambda|=\sum_{1\leq
i\leq n}|\lambda_i|$. We define
\begin{equation}
\lambda^*=(-\lambda_n,\ldots,-\lambda_1).
\end{equation}
Clearly, $\lambda^*\in\mathbb{Z}_+^n$, and $\lambda^{**}=\lambda$.
We put ${\bf 0}_n=(0,\ldots,0)\in \mathbb{Z}_+^n$.

Set $[-n]= \{\, -n<\cdots < -1 \}$ with $[-n]_0=[-n]$. Given
$\lambda=(\lambda_1,\cdots,\lambda_n)\in \mathbb{Z}_+^n$, we may
identify $\lambda$ with a {\it generalized Young diagram} in the
following way. First, we fix a vertical line. Then for each
$\lambda_{k}$, we place $|\lambda_k|$ nodes in the $k$th row in a
left-justified (resp. right-justified) way with respect to the
vertical line if $\lambda_{k}\geq 0$ (resp. $\lambda_k\leq 0$). For
example,\vskip 3mm
\begin{center}
 $\lambda=(3,2,0,-1,-2) \longleftrightarrow$
\begin{tabular}{cc|ccc}
   &   & $\bullet$ & $\bullet$ & $\bullet$  \\
   &   & $\bullet$ & $\bullet$ &   \\
  - & - & - & -& - \\
   & $\bullet$ &   &   &  \\
 $\bullet$ & $\bullet$ &   & &   \\
\multicolumn{5}{c}{\small -2  \ -1 \ \  1 \ \ 2 \ \ 3  }
\end{tabular}\ \ \ \ ,
\end{center}
where $-$ denotes the empty row. We enumerate the columns of a
diagram as in the above figure. We put
\begin{equation}
\begin{split}
\lambda^+&=(\max\{\lambda_1,0\},\cdots,\max\{\lambda_n,0\}), \\
\lambda^-&=(\max\{-\lambda_n,0\},\cdots,\max\{-\lambda_1,0\}),
\end{split}
\end{equation}
where $\lambda^{\pm}$ are understood as partitions (or Young
diagrams).

\begin{df}[\cite{St}]\label{rationalSST}{\rm
Let $T$ be a tableau obtained by filling a generalized Young diagram
$\lambda$ of level $n$ with entries in $[n] \cup [-n]$. We call $T$
a {\it rational $[n]$-semistandard tableau of shape $\lambda$} if
\begin{itemize}
\item[(1)] the entries in the columns indexed by positive (resp. negative)
numbers belong to $[n]$ (resp. $[-n]$),
\item[(2)] the entries in each row (resp. column) are weakly (resp. strictly)
increasing from left to right (resp. from top to bottom),
\item[(3)] if $b_1<\cdots<b_s$ (resp. $-b'_1<\cdots<-b'_t$) are the entries
in the $1$st (resp. $-1$st) column ($s+t\leq n$), then
$$b''_i\leq b_i,$$ for $1\leq i\leq s$, where $\{ b''_1<\cdots<b''_{n-t} \}=
[n]\setminus \{ b'_1,\cdots,b'_t \}$.
\end{itemize}}
\end{df}
We denote by $SST_{[n]}(\lambda)$ the set of all rational
$[n]$-semistandard tableaux of shape $\lambda\in\mathbb{Z}_+^n$.

\begin{ex}{\rm For $\lambda=(3,2,0,-1,-2)$, we have

\begin{center}
\begin{tabular}{cc|ccc}
   &   & $2$ & $3$ & $5$ \\
   &   & $4$ & $4$ &  \\
  - & -&- & -& -\\
   & $-5$ &   &   &  \\
 $-3$ & $-3$ &   &   &
\end{tabular}$\in SST_{[5]}(\lambda)$.
\end{center}}
\end{ex}

Let us explain the relation between the rational $[n]$-semistandard
tableaux and ordinary $[n]$-semistandard tableaux. Let $T$ be a
rational $[n]$-semistandard tableau of shape $(0^{n-t},(-1)^{t})$
($0\leq t\leq n$) with the entries $-b_1<\cdots<-b_t$. We define
$\sigma(T)$ to be the $[n]$-semistandard tableau of shape
$(1^{n-t},0^{t})$ with the entries $b'_1<\cdots<b'_{n-t}$, where $\{
b'_1<\cdots<b'_{n-t}\}= [n]\setminus \{b_1<\cdots<b_t\}$. If $t=n$,
then we define $\sigma(T)$ to be the empty tableau.

Generally, for $\lambda\in\mathbb{Z}_+^n$ and $T\in
SST_{[n]}(\lambda)$, we define $\sigma(T)$ to be the tableau
obtained by applying $\sigma$ to the $-1$st column of $T$. For
example, when $n=5$, we have \vskip 3mm
\begin{center}
$\sigma\left(\begin{array}{cc|ccc}
   &   & 2 & 3 & 5 \\
   &   & 4 & 4 &  \\
  \text{-} & \text{-}& \text{-}&\text{-} & \text{-}\\
   & {\bf -5} &   &   &  \\
 -4 & {\bf -2} &   &   &
\end{array}\right)$\ \ =
\begin{tabular}{c|cccc}
   & {\bf $\bf 1$} & $2$ & $3$ & $5$ \\
   & {\bf $\bf 3$} & $4$ & $4$ &  \\
   & {\bf $\bf 4$} & & & \\
  - & - &  - &  - & - \\
 $-4$ &  &   &   &
\end{tabular}.
\end{center}\vskip 3mm
By Definition \ref{rationalSST}, it is straightforward to check that
$\sigma(T)\in SST_{[n]}(\lambda+(1^n))$, where we define
$\mu+\nu=(\mu_k+\nu_k)_{k\geq 1}$ for $\mu=(\mu_k)_{k\geq 1}$ and
$\nu=(\nu_k)_{k\geq 1}$ in $\mathbb{Z}^n_+$.

\begin{lem}\label{shift of SST}
For $\lambda\in\mathbb{Z}^n_+$,  the map $\sigma :
SST_{[n]}(\lambda) \rightarrow SST_{[n]}(\lambda+(1^n))$ is a
bijection.\qed
\end{lem}

Next, let us introduce the notion of a {\it rectangular complement}
of an $[n]$-semistandard tableau. Fix $n\geq 1$. For $\lambda\in
\cP_n$ and a tableau $T\in SST_{[n]}(\lambda)$, we define
\begin{equation}
\delta^n_k(T)=(\sigma^{-k}(T))^{\pi},
\end{equation}
for $k\geq \lambda_1$. We put
\begin{equation}\label{delta nk lambda}
\delta^n_k(\lambda)=(k-\lambda_n,\ldots,k-\lambda_1).
\end{equation}
Then $\delta^n_k(\lambda)$ is the shape of $\delta^n_k(T)$. By
definition, we have $\delta^n_k(T)\in
SST_{[-n]^{\pi}}(\delta^n_k(\lambda))$. Identifying $[-n]^{\pi}$
with $[n]$ ($-i$ with $i$ for each $i$), we may view $\delta^n_k(T)$
as an element in $SST_{[n]}(\delta^n_k(\lambda))$.

\begin{ex}{\rm If $\lambda=(4,3,1)\in\cP_4$, and $T=$
\begin{tabular}{cccc}
   1 & $2$ & $2$ & $3$ \\
   3 & $4$ & $4$ &  \\
   4 & & &
\end{tabular}
, then
$$\sigma^{-5}(T)=
\begin{tabular}{ccccc}
    &  &  &   & -4 \\
    &  &  & -4 & -3\\
    & -3 & -3 & -2 & -2 \\
   -2 & -1 & -1 & -1 & -1
\end{tabular}\ \ \ \ \text{and} \ \ \ \
\delta^4_5(T)=[\sigma^{-5}(T)]^{\pi}=
\begin{tabular}{ccccc}
  1 & 1 & 1 & 1 & 2 \\
  2 & 2 & 3 & 3 &  \\
  3 & 4 & & & \\
  4 & & & &
\end{tabular}.$$
}
\end{ex}\vskip 3mm

Since $\delta^n_k\circ \delta^n_k (T)=T$ for all $T\in
SST_{[n]}(\lambda)$, it follows that
\begin{lem}\label{delta nk} For $\lambda\in\cP_n$ and $k\geq \lambda_1$, the map
$$\delta^n_k : SST_{[n]}(\lambda)\longrightarrow
SST_{[n]}(\delta^n_k(\lambda))$$ is a bijection. \qed
\end{lem}
\vskip 3mm

The row (resp. column) insertion of $[n]$-semistandard tableaux
corresponds to the column (resp. row) insertion  of their
rectangular complements as we can see in the following theorem due
to Stroomer \cite{Str}. This fact will be of important use in the
proof of our main results.

\begin{thm}[\cite{Str} Theorem 5.11]\label{Stroomer}
Assume that $\mu$ and $\nu \in \cP_n$ are given. For $T_1\in
SST_{[n]}(\mu)$ and $T_2\in SST_{[n]}(\nu)$, we have
$$\delta_{p+q}^n(T_2\rightarrow T_1)=(\delta^n_{p}(T_1)\leftarrow
\delta^n_{q}(T_2)),$$ for $p\geq \mu_1$ and $q\geq \nu_1$. \qed
\end{thm}\vskip 3mm

For $n\geq 1$, let ${\bf x}_{[n]}=\{\,x_1,\ldots,x_n\,\}$ be the set
of $n$ variables. For $\lambda\in\mathbb{Z}_+^n$, we define the
character of $SST_{[n]}(\lambda)$ to be
\begin{equation}
s_{\lambda}({\bf x}_{[n]})=\sum_{T\in SST_{[n]}(\lambda)}{\bf
x}_{[n]}^{T},
\end{equation}
where ${\bf x}_{[n]}^{T}=\prod_{k\in [n]}x_k^{m_k-m_{-k}}$ and $m_k$
is the number of occurrences of $k$ in $T$ for $k\in [n]\cup [-n]$.
Note that $s_{{\bf 0}_n}({\bf x}_{[n]})=1$. We call
$s_{\lambda}({\bf x}_{[n]})$ the {\it rational Schur polynomial
corresponding to $\lambda$}. Note that $s_{\lambda}({\bf x}_{[n]})$
is the ordinary Schur polynomial if $\lambda\in\cP_n$, and
$(x_1\cdots x_n)^k s_{\lambda}({\bf x}_{[n]})=s_{\lambda+(k^n)}({\bf
x}_{[n]})$ for $k\geq 1$.

By Theorem \ref{LR rule} and the linear independence of Schur
polynomials, it follows that for $\mu,\nu\in\mathbb{Z}_+^n$,
\begin{equation}\label{LR rule for rational Schur}
s_{\mu}({\bf x}_{[n]})s_{\nu}({\bf
x}_{[n]})=\sum_{\lambda\in\mathbb{Z}_+^n}c^{\lambda}_{\mu\,\nu}s_{\lambda}({\bf
x}_{[n]}),
\end{equation}
where
$c^{\lambda}_{\mu\,\nu}=N^{\lambda+((p+q)^n)}_{\mu+(p^n)\,\nu+(q^n)}$
for all sufficiently large $p,q>0$. By definition, we also have
$c^{\lambda}_{\mu\,\nu}=c^{\lambda}_{\nu\,\mu}$. For later use, let
us characterize $c^{\lambda}_{\mu\,\nu}$ more explicitly in terms of
Littlewood-Richardson tableaux.

\begin{lem} Suppose that $\lambda,\mu$ and
$\nu\in\cP_n$ are given. For $p,q\geq 0$, there exists a bijection
$$\pi^n_{p,q} : LR^{\lambda}_{\mu\,\nu}\longrightarrow LR^{\lambda+((p+q)^n)}_{\mu+(p^n)\,\nu+(q^n)}.$$
\end{lem}
\pf Suppose that $LR^{\lambda}_{\mu\,\nu}$ is non-empty and $Q\in
LR^{\lambda}_{\mu\,\nu}$ is given. Clearly, we may view $Q$ as an
element in $LR^{\lambda+(k^{n})}_{\mu+(k^n)\,\nu}$ for $k\geq 0$.
This implies the bijection
$$\pi^n_{k,0} : LR^{\lambda}_{\mu\,\nu}
\longrightarrow LR^{\lambda+(k^{n})}_{\mu+(k^n)\,\nu}.$$ Then
$\pi^n_{p,q}$ is given by $\pi^n_{p,q}=\theta\circ
\pi^n_{q,0}\circ\theta \circ \pi^n_{p,0}$, where $\theta$ denotes
the map given in Corollary \ref{skew LR} (2). \qed\vskip 3mm

Suppose that $\lambda,\mu$ and $\nu\in\mathbb{Z}_+^n$ are given. Let
$p$ and $q$ be the smallest non-negative integers such that
$\mu+(p^n),\nu+(q^n),\lambda+((p+q)^n)\in\cP_n$. Then we define
\begin{equation}\label{LR tableaux for rational}
{\bf LR}^{\lambda/\mu}_{\nu}=\{\,[Q]\,|\,Q\in
LR^{\lambda+((p+q)^{n})}_{\mu+(p^n)\,\nu+(q^n)} \,\},
\end{equation}
where $[Q]=\{\,\pi^n_{s,t}(Q)\,|\,s,t\geq 0\,\}$.
\begin{lem}
Under the above hypothesis, $|{\bf
LR}^{\lambda/\mu}_{\nu}|=c^{\lambda}_{\mu\,\nu}$
\end{lem}
\pf It follows directly from \eqref{LR rule for rational
Schur}.\qed\vskip 3mm

Next, if we put ${\bf x}_{[m+n]}={\bf x}_{[m]}\sqcup {\bf y}_{[n]}$
for $m,n>0$, where ${\bf y}_{[n]}=\{\,y_i=x_{m+i}\,|\,i\in [n]\,\}$,
then by Corollary \ref{skew LR} (1), we have for
$\lambda\in\mathbb{Z}_+^{m+n}$
\begin{equation}\label{LR' rule for rational Schur}
s_{\lambda}({\bf x}_{[m+n]})=s_{\lambda}({\bf x}_{[m]},{\bf
y}_{[n]}) =\sum_{\mu,\nu}\hat{c}^{\lambda}_{\mu\,\nu}s_{\mu}({\bf
x}_{[m]})s_{\nu}({\bf y}_{[n]}),
\end{equation}
where
$\hat{c}^{\lambda}_{\mu\,\nu}=N^{\lambda+(p^{m+n})}_{\mu+(p^m)\,\nu+(p^n)}$
for all sufficiently large $p>0$. We may also characterize
$\hat{c}^{\lambda}_{\mu\,\nu}$ in terms of Littlewood-Richardson
tableaux.

\begin{lem}\label{LR coeff for SAB} For $\lambda\in\cP_{m+n},\mu\in\cP_m$ and $\nu\in \cP_n$, there exists a bijection
$$\pi^{m,n}_{\ell} : LR^{(\lambda+(\ell^{m+n}))'}_{(\mu+(\ell^m))'\,(\nu+(\ell^n))'}
\longrightarrow LR^{\lambda'}_{\mu'\,\nu'},$$ where $\ell\geq 0$.
\end{lem}
\pf Suppose that
$LR^{(\lambda+(\ell^{m+n}))'}_{(\mu+(\ell^m))'\,(\nu+(\ell^n))'}$ is
non-empty and $Q\in
LR^{(\lambda+(\ell^{m+n}))'}_{(\mu+(\ell^m))'\,(\nu+(\ell^n))'}$ is
given. Note that the $i$th row of $Q$ is filled with $i$ for $1\leq
i\leq \ell$ (in fact, the first $\ell$ rows of $Q$ is of the form
$H^{(n^{\ell})}$), and the other entries in $Q$ are greater than
$\ell$ since the content of $Q$ is $(\nu+(\ell^n))'$. Now, we define
$\pi^{m,n}_{\ell}(Q)$ to be the tableau obtained from $Q$ by
\begin{itemize}
\item[(1)] removing the first $\ell$ rows of $Q$,

\item[(2)] replacing each entry $i$ in the remaining part of $Q$ by
$i-\ell$.
\end{itemize}
Then it is not difficult to see that $\pi^{m,n}_{\ell}(Q)\in
LR^{\lambda'}_{\mu'\,\nu'}$. Since the first $\ell$ rows in any
$Q\in
LR^{(\lambda+(\ell^{m+n}))'}_{(\mu+(\ell^m))'\,(\nu+(\ell^n))'}$
is of the form $H^{(n^{\ell})}$, any $Q'\in
LR^{\lambda'}_{\mu'\,\nu'}$ is given by $\pi^{m,n}_{\ell}(Q)$ for
a unique $Q\in
LR^{(\lambda+(\ell^{m+n}))'}_{(\mu+(\ell^m))'\,(\nu+(\ell^n))'}$.
Hence $\pi^{m,n}_{\ell}$ is a bijection. \qed\vskip 3mm

Suppose that $\mu\in\mathbb{Z}_+^m$, $\nu\in\mathbb{Z}_+^n$, and
$\lambda\in\mathbb{Z}_+^{m+n}$ are given. Let $d$ be the smallest
non-negative integer such that $\lambda+(d^{m+n}), \mu+(d^m),
\nu+(d^n)\in \cP$. Then we define
\begin{equation}\label{skew LR tableaux for rational}
{\bf LR}^{\lambda}_{\mu\,\nu}=\{\,[Q]\,|\,Q\in
LR^{(\lambda+(d^{m+n}))'}_{(\mu+(d^m))'\,(\nu+(d^n))'}\,\},
\end{equation}
where $[Q]$ is the set of all the bijective images
$(\pi^{m,n}_{\ell})^{-1}(Q)$ ($\ell\geq 0$). From \eqref{LR' rule
for rational Schur}, it follows that
\begin{lem}
Under the above hypothesis,  $|{\bf
LR}^{\lambda}_{\mu\,\nu}|=\hat{c}^{\lambda}_{\mu\,\nu}$.
\end{lem}
\pf It follows from \eqref{LR' rule for rational Schur} and
Corollary \ref{conjugate LR}.\qed\vskip 5mm

\subsection{$\A/\B$-semistandard tableaux}
Now, let us introduce a certain class of bitableaux, which is our
main object in this paper.

\begin{df}\label{ABsst}{\rm
Suppose that $\lambda\in\mathbb{Z}_+^n$ is given. An {\it
$\A/\B$-semistandard tableau of shape $\lambda$} is a pair of
tableaux $(T^+,T^-)$ such that
\begin{equation*}
T^+\in SST_{\A}((\lambda+(d^n))/\mu), \ \ \ \ T^- \in
SST_{\B}((d^{n})/\mu),
\end{equation*}
for some integer $d\geq 0$ and $\mu\in\cP_n$ satisfying
\begin{itemize}
\item[(1)] $\lambda+(d^n)\in\cP_n$,

\item[(2)] $\mu\subset (d^n)$, and $\mu\subset\lambda+(d^n)$.
\end{itemize}}
\end{df}\vskip 3mm
We denote by $SST_{\A/\B}(\lambda)$ the set of all
$\A/\B$-semistandard tableaux of shape $\lambda$.  For $(T^+,T^-)\in
SST_{\A/\B}(\lambda)$, we say that $(T^+,T^-)$ is of {\it level
$n$}. We set
\begin{equation}
\cP_{\A/\B}=\bigsqcup_{n\geq 1}\cP_{\A/\B,n}
\end{equation}
where
$\cP_{\A/\B,n}=\{\,\lambda\in\mathbb{Z}_+^n\,|\,SST_{\A/\B}(\lambda)\neq
\emptyset\,\}$

\begin{ex}{\rm
Suppose that $\A=\B=\mathbb{N}$. Consider
\begin{equation*}
(T^+,T^-)= \left(
\begin{array}{ccc|ccc}
   & & 1& 1& 2& 2\\
   & 1& 2& 2& 4& \\
  2 & 3& 3& & & \\
  4 & & & & & \\
\end{array}
,
\begin{array}{ccc|ccc}
   &  &  1& & & \\
   &  1&  2& & & \\
  2 & 2 & 4&  & & \\
  3 & 3 & 5&  & & \\
\end{array}
\right)\in SST_{\A/\B}((3,2,0,-2)),
\end{equation*}
where the vertical lines in $T^+$ and $T^-$ correspond to the one in
the generalized partition $(3,2,0,-2)$. Note that
\begin{equation*}
\begin{split}
{\rm sh}(T^+)&=\left((3,2,0,-2)+(3^4)\right)/(2,1,0,0), \\
{\rm sh}(T^-)&=(3^4)/(2,1,0,0).
\end{split}
\end{equation*}
  }
\end{ex}\vskip 3mm

\begin{ex}{\rm
Let $\lambda={\bf 0}_n\in\mathbb{Z}_+^n$. Then for $(T^+,T^-)\in
SST_{\A/\B}({\bf 0}_n)$, $${\rm sh}(T^+)={\rm sh}(T^-)=(d^n)/\mu$$
for some $d\geq 0$ and $\mu\in\cP_n$. If we identify $(T^+,T^-)$
with $((T^+)^{\pi},(T^-)^{\pi})$, then we have
\begin{equation}\label{SST for 0}
SST_{\A/\B}({\bf 0}_n)=\bigsqcup_{\lambda\in
\cP_n}SST_{\A}(\lambda)\times SST_{\B}(\lambda),
\end{equation}
by Lemma \ref{reverse ordering}.}
\end{ex}\vskip 3mm

%\begin{rem}{\rm
%We assume that for $(T^+,T^-)\in SST_{\A/\B}(\lambda)$, either $T^+$
%or $T^-$ may be an empty tableau, but they are not both empty. For
%example, for $\lambda\in\cP_n$, a tableau $T\in SST_{\A}(\lambda)$
%can be identified with $(T^+,\emptyset)\in SST_{\A/\B}(\lambda)$.

%(2) The notion of $\A/\B$-semistandard tableaux might be viewed as a
%generalization of rational $[n]$-semistandard tableaux in
%\ref{rational tableaux} in the following sense: for
%$\lambda\in\mathbb{Z}_+^n$, a tableau $T\in SST_{[n]}(\lambda)$ may
%be identified with the pair $(T^+,T^-)\in SST_{[n]/[-n]}(\lambda)$,
%where $T^+$ (resp. $T^-$) is the part of $T$ with positive (resp.
%negative) column indices. Note that
%\begin{equation*}
%\begin{split}
%& {\rm sh}(T^+)=(\lambda_1,\ldots,\lambda_r), \\
%& {\rm sh}(T^-)^{\pi}=(-\lambda_n,\ldots,-\lambda_{r+1}),
%\end{split}
%\end{equation*}
%if $\lambda=(\lambda_1\geq\cdots \geq\lambda_r\geq
%0>\lambda_{r+1}\geq\cdots\geq\lambda_n$). However,
%$SST_{[n]}(\lambda)$ is a proper subset of
%$SST_{[n]/[-n]}(\lambda)$}.
%}
%\end{rem}

The decomposition \eqref{SST for 0} can be viewed as a branching
rule, and it can be generalized as follows.
\begin{prop}\label{branching for AB tableaux}
For  $\lambda\in \cP_{\A/\B, n}$, there exists a bijection between
$SST_{\A/\B}(\lambda)$ and
$$\bigsqcup_{\mu,\nu\in\cP_n}{\bf LR}^{\lambda/\mu}_{\nu^*}\times
SST_{\A}(\mu)\times SST_{\B}(\nu).
$$
\end{prop}
\pf By Definition \ref{ABsst}, we have
\begin{equation*}
\begin{split}
&SST_{\A/\B}(\lambda) \\ & \stackrel{1-1}{\longleftrightarrow}
\bigsqcup_{\substack{d\geq 0
\\ \eta\in \cP_n}}SST_{\A}((\lambda+(d^n))/\eta)\times
SST_{\B}((d^n)/\eta) \\
& \stackrel{1-1}{\longleftrightarrow} \bigsqcup_{\substack{d\geq 0
\\ \eta,\mu \in \cP_n}} LR^{\lambda+(d^n)}_{\mu\,\eta} \times
SST_{\A}(\mu) \times
SST_{\B}((d^n)/\eta) \ \ \ \ \ \text{by Corollary \ref{skew LR}} \\
& \stackrel{1-1}{\longleftrightarrow} \bigsqcup_{\substack{d\geq 0
\\ \eta,\mu \in \cP_n}} {\bf LR}^{\lambda/\mu}_{\eta-(d^n)} \times
SST_{\A}(\mu) \times SST_{\B}((d^n)/\eta)  \ \ \ \ \ \text{by
\eqref{LR tableaux for rational}},
\end{split}
\end{equation*}
where the union is taken over $d\geq 0$ and $\eta\in\cP_n$, which
satisfy the conditions (1) and (2) in Definition \ref{ABsst}. For
each $d\geq 0$ and $\eta\in\cP_n$ above, there exists a unique
$\nu\in\cP_n$ such that $(d^n)/\eta=\nu^{\pi}$, and
$\eta-(d^n)=\nu^{\ast}$. Hence, we have
$$SST_{\B}((d^n)/\eta) \stackrel{1-1}{\longleftrightarrow} SST_{\B^{\pi}}(\nu)
\stackrel{1-1}{\longleftrightarrow} SST_{\B}(\nu),$$ where the first
correspondence is given by $\pi$, and the second one is given by
Proposition \ref{reverse ordering}. This completes the proof.
\qed\vskip 5mm

\subsection{Characters}

From now on, let ${\bf x}_{\A}=\{\,x_a\,|\,a\in\A\,\}$ be the set of
variables indexed by $\A$. Let $P_{\A}=\bigoplus_{a\in
\A}\mathbb{Z}\epsilon_a$ be the free abelian group with the basis
$\{\,\epsilon_a\,|\,a\in \A\,\}$. For $\lambda\in \cP_{\A}$ and
$T\in SST_{\A}(\lambda)$, we define the {\it $\A$-weight of $T$} by
${\rm wt}_{\A}(T)=\sum_{a\in\A}m_a\epsilon_a\in P_{\A}$ where $m_a$
is the number of occurrences of $a$ in $T$. Put ${\bf
x}_{\A}^{T}=\prod_{a\in\A}x_a^{m_a}$. Now, we define the {\it
character of $SST_{\A}(\lambda)$} to be
\begin{equation}
S_{\lambda}({\bf x}_{\A})=\sum_{T\in SST_{\A}(\lambda)}{\bf
x}_{\A}^{T}\in \mathbb{Z}[[{\bf x}_{\A}]].
\end{equation}
For simplicity, let us often write
$S^{\A}_{\lambda}=S_{\lambda}({\bf x}_{\A})$. We assume that
$S^{\A}_{\lambda}=0$ unless $\lambda\in\cP_{\A}$. The character of
$SST_{\A}(\lambda/\mu)$ of a skew Young diagram is defined
similarly.

\begin{ex}\label{ex of characters}{\rm
(1) If $\A=[n]$ and $\lambda\in\cP_n$, then $S_{\lambda}({\bf
x}_{[n]})=s_{\lambda}({\bf x}_{[n]})$ is the {\it Schur polynomial
corresponding to $\lambda$} (cf.\cite{Mac}).

(2) Suppose that $\A=[m]\ast [n]'$. To distinguish $[m]$ and $[n]'$
as sets, let us write $[n]'=\{1'<2'<\cdots<n'\}$. Then
$S_{\lambda}({\bf x}_{\A})$ is the {\it $(m,n)$-hook Schur
polynomial corresponding to $\lambda$}, which is the character of an
irreducible representation of the general linear Lie superalgebra
$\frak{gl}_{m|n}$ (see \cite{BR}). Note that $S_{\lambda}({\bf
x}_{\A})\neq 0$ if and only if $\lambda_{m+1}\leq n$, that is,
$\lambda$ is an {\it $(m,n)$-hook partition}. }
\end{ex}\vskip 3mm

\begin{lem}\label{switching character}
For $\lambda\in \cP_{\A}$,
$S^{\A}_{\lambda}=S^{\A_0\ast\A_1}_{\lambda}=S^{\A^{\pi}}_{\lambda}$.
\end{lem}
\pf It follows directly from Proposition \ref{shuffling} and
\ref{reverse ordering}. \qed\vskip 3mm

\begin{lem}\label{linear indep of SA} The set $\{\,S^{\A}_{\lambda}\,|\,\lambda\in \cP_{\A}\,\}$
is linearly independent over $\mathbb{Z}$.
\end{lem}
\pf Suppose that
\begin{equation*}
\sum_{1\leq i\leq m}a_{\lambda^{(i)}}S^{\A}_{\lambda^{(i)}}=0,
\end{equation*}
where $\lambda^{(i)}\in \cP_{\A}$ and
$a_{\lambda^{(i)}}\in\mathbb{Z}$ for $1\leq i\leq m$.

Choose a finite subset $\tilde{\A}\subset \A$ such that
$S^{\tilde{\A}}_{\lambda^{(i)}}\neq 0$ for all $1\leq i\leq m$. If
we put $x_a=0$ for all $a\not\in\tilde{\A}$ in
$S^{\A}_{\lambda^{(i)}}$, then we get
\begin{equation*}
\sum_{1\leq i\leq
n}a_{\lambda^{(i)}}S^{\tilde{\A}}_{\lambda^{(i)}}=0.
\end{equation*}
By Lemma \ref{switching character}, $S^{\tilde{\A}}_{\lambda^{(i)}}$
are hook Schur polynomials (see Example \ref{ex of characters}).
Then the linear independence follows from that of ordinary hook
Schur polynomials (see \cite{BR}). \qed\vskip 3mm

Since the bijections in Theorem \ref{LR rule} and Corollary
\ref{skew LR} (1) preserve the weights of $\A$-semistandard
tableaux, it follows that
\begin{cor}\label{branching LR}\mbox{}
\begin{itemize}
\item[(1)] For $\mu,\nu \in \cP_{\A}$, we have
$S^{\A}_{\mu}S^{\A}_{\nu}=\sum_{\lambda\in\cP_{\A}}N^{\lambda}_{\mu\,\nu}S^{\A}_{\lambda}$.

\item[(2)] For a skew Young diagram $\lambda/\mu$, we have
$S^{\A}_{\lambda/\mu}=\sum_{\nu\in\cP_{\A}}N^{\lambda}_{\mu\,\nu}S^{\A}_{\nu}$.
\end{itemize}\qed
\end{cor}
In particular, for $\lambda\in \cP_{\A/\B}$, we have
\begin{equation*}
S^{\A\ast\B}_{\lambda}=\sum_{\mu\in
\cP}S^{\A}_{\mu}S^{\B}_{\lambda/\mu}
=\sum_{\mu\in\cP_{\A},\nu\in\cP_{\B}}N^{\lambda}_{\mu\,\nu}S^{\A}_{\mu}S^{\B}_{\nu}.
\end{equation*}\vskip 3mm

Now, for $\lambda\in\mathbb{Z}_+^n$, we define the {\it character of
$SST_{\A/\B}(\lambda)$} to be
\begin{equation}
S_{\lambda}({\bf x}_{\A},{\bf x}_{\B})=\sum_{(T^+,T^-)\in
SST_{\A/\B}(\lambda)}{\bf x}_{\A}^{T^+}({\bf x}_{\B}^{T^-})^{-1}\in
\mathbb{Z}[[{\bf x}_{\A},{\bf x}_{\B}^{-1}]],
\end{equation}
where ${\bf x}^{-1}_{\B}=\{\,x_b^{-1}\,|\,b\in\B\,\}$. For
simplicity, let us write $S^{\A/\B}_{\lambda}=S_{\lambda}({\bf
x}_{\A},{\bf x}_{\B})$. We assume that $S^{\A/\B}_{\lambda}=0$
unless $\lambda\in \cP_{\A/\B}$. It is easy to check that
$S^{\A/\B}_{\lambda}$ is a well-defined element in $\mathbb{Z}[[{\bf
x}_{\A},{\bf x}_{\B}^{-1}]]$.

Note that each non-zero monomial in $S^{\A/\B}_{\lambda}$ is of
homogeneous degree $\langle\lambda\rangle=\sum_{i=1}^n\lambda_i$,
and the degree of each monomial in ${\bf x}_{\A}$ (resp. ${\bf
x}_{\B}^{-1}$) is at least $|\lambda^+|$ (resp. $|\lambda^-|$). For
$(T^+,T^-)\in SST_{\A/\B}(\lambda)$, we also define the {\it
$\A/\B$-weight of $T$} by ${\rm wt}_{\A/\B}(T)={\rm
wt}_{\A}(T^+)-{\rm wt}_{\B}(T^-)\in P_{\A}\oplus P_{\B}$.

\begin{rem}\label{formal sum}{\rm
If we identify $\mathbb{Z}[[{\bf x}_{\A},{\bf x}_{\B}^{-1}]]$ with
$\mathbb{Z}[[{\bf x}_{\A}]]\otimes\mathbb{Z}[[{\bf x}_{\B}^{-1}]]$,
then the set $$S_{\A,\B}=\{\,S_{\lambda}({\bf x}_{\A})S_{\mu}({\bf
x}_{\B}^{-1})\,|\,\lambda\in \cP_{\A},\mu\in\cP_{\B}\,\}$$ is
linearly independent over $\mathbb{Z}$. Consider
$$\sum_{\lambda,\mu}c_{\lambda\,\mu}S_{\lambda}({\bf
x}_{\A})S_{\mu}({\bf x}_{\B}^{-1}),$$ for
$c_{\lambda\,\mu}\in\mathbb{Z}$, which is not necessarily a finite
sum. Then we can check that it is a well-defined element in
$\mathbb{Z}[[{\bf x}_{\A},{\bf x}_{\B}^{-1}]]$ since each monomial
in ${\bf x}_{\A}$ and ${\bf x}_{\B}$ occurs only in finitely many
$\lambda$ and $\mu$'s. From the linear independence of $S_{\A,\B}$,
$c_{\lambda\,\mu}$ is uniquely determined for all $\lambda,\mu$. }
\end{rem}

Now, we can express $S_{\lambda}({\bf x}_{\A},{\bf x}_{\B})$ as a
(possibly infinite) linear combination of $S_{\lambda}({\bf
x}_{\A})S_{\mu}({\bf x}_{\B}^{-1})$'s as follows.
\begin{prop}\label{character for SAB} For  $\lambda\in\cP_{\A/\B,n}$, we have
\begin{equation*}
S_{\lambda}({\bf x}_{\A},{\bf x}_{\B})=\sum_{\mu,\nu\in
\cP_n}c^{\lambda}_{\mu\,\nu^*}S_{\mu}({\bf x}_{\A})S_{\nu}({\bf
x}^{-1}_{\B}).
\end{equation*}
\end{prop}
\pf It follows from Proposition \ref{branching for AB tableaux}.
\qed

\begin{ex}{\rm  When $\lambda={\bf 0}_n$, we have
\begin{equation*}
S_{{\bf 0}_n}({\bf x}_{\A},{\bf
x}_{\B})=\sum_{\lambda\in\cP_n}S_{\lambda}({\bf
x}_{\A})S_{\lambda}({\bf x}_{\B}^{-1}).
\end{equation*}}
\end{ex}

\begin{prop}\label{linear indep of SAB}
The set $\{\,S^{\A/\B}_{\lambda}\,|\,\lambda\in  \cP_{\A/\B} \,\}$
is linearly independent over $\mathbb{Z}$.
\end{prop}
\pf Suppose that
\begin{equation*}
\sum_{1\leq i\leq m}a_{\lambda^{(i)}}S^{\A/\B}_{\lambda^{(i)}}=0,
\end{equation*}
where $\lambda^{(i)}\in \cP_{\A/\B}$ and
$a_{\lambda^{(i)}}\in\mathbb{Z}$ for $1\leq i\leq m$. Suppose that
the level of $\lambda^{(i)}$ is $n_i$ for $1\leq i\leq m$. We will
use induction on $m$ to show that $a_{\lambda^{(i)}}=0$ for $1\leq
i\leq m$. Clearly, it is true when $m=1$.

Let $n=\max\{n_1,\ldots,n_m\}$, and let $I=\{\,i\,|\,n_i=n\,\}$.
Choose a positive integer $d$ such that
$\lambda^{(i)}+(d^{n})\in\cP$ for all $i\in I$. Then if we write
$S^{\A/\B}_{\lambda^{(i)}}$ $(i\in I)$ as a linear combination of
$S_{\lambda}({\bf x}_{\A})S_{\mu}({\bf x}_{\B}^{-1})$'s, then the
coefficient of $S_{(d^{n})}({\bf x}_{\B}^{-1})$ in
$S^{\A/\B}_{\lambda^{(i)}}$ is $S_{\lambda^{(i)}+(d^{n})}({\bf
x}_{\A})$ for $i\in I$ (cf. Proposition \ref{branching for AB
tableaux} and \ref{character for SAB}). Since $S_{(d^{n})}({\bf
x}_{\B}^{-1})$ occurs only in the expansion of
$S^{\A/\B}_{\lambda^{(i)}}$ for $i\in I$, we have
\begin{equation*}
\sum_{i\in I}a_{\lambda^{(i)}}S_{\lambda^{(i)}+(d^n)}({\bf
x}_{\A})=0,
\end{equation*}
from the linear independence of $S_{\A,\B}$ (see also Remark
\ref{formal sum}). Since $\lambda^{(i)}+(d^n)$ are mutually
different for $i\in I$, we have $a_{\lambda^{(i)}}=0$ for $i\in
I$, and hence by induction hypothesis, $a_{\lambda^{(i)}}=0$ for
all $1\leq i\leq m$. \qed\vskip 3mm

\begin{cor}\label{linear indep of SAB'}
For $n\geq 1$, consider
$\sum_{\lambda\in\cP_{\A/\B,n}}a_{\lambda}S^{\A/\B}_{\lambda}$,
where $a_{\lambda}\in\mathbb{Z}$, which is not necessarily a finite
sum. Then it is a well-defined element in $\mathbb{Z}[[{\bf
x}_{\A},{\bf x}_{\B}^{-1}]]$, and the coefficient $a_{\lambda}$ is
uniquely determined for $\lambda\in\cP_{\A/\B,n}$.
\end{cor}
\pf Note that for $\mu,\nu\in\cP_n$, $S_{\mu}({\bf
x}_{\A})S_{\nu}({\bf x}_{\B}^{-1})$ occurs in the expansion of
$S^{\A/\B}_{\lambda}$ only if
\begin{equation}\label{formal sum condition}
|\mu|\geq |\lambda^+|\ \ \ \text{and} \ \ \ |\nu|\geq |\lambda^-| \
\ \ \ \ (\ \ \text{or $|\mu|+|\nu|\geq |\lambda|$}\ \ ),
\end{equation}
Since the level of $\lambda$ is fixed, there are only  finitely many
$\lambda$'s satisfying \eqref{formal sum condition}, and hence the
coefficient of $S_{\mu}({\bf x}_{\A})S_{\nu}({\bf x}_{\B}^{-1})$ in
$\sum_{\lambda\in\cP_{\A/\B,n}}a_{\lambda}S^{\A/\B}_{\lambda}$ is a
well-defined integer. This implies that
$\sum_{\lambda\in\cP_{\A/\B,n}}a_{\lambda}S^{\A/\B}_{\lambda}$ is a
well-defined element in $\mathbb{Z}[[{\bf x}_{\A},{\bf
x}_{\B}^{-1}]]$ (see Remark \ref{formal sum}).

Next, suppose that
$\sum_{\lambda\in\cP_{\A/\B,n}}a_{\lambda}S^{\A/\B}_{\lambda}=0$.
For $d>0$, the coefficient of $S_{(d^n)}({\bf x}_{\B}^{-1})$ in
$\sum_{\lambda\in\cP_{\A/\B,n}}a_{\lambda}S^{\A/\B}_{\lambda}$ is
given by
$$\sum_{\lambda+(d^n)\in\cP_n}a_{\lambda}S_{\lambda+(d^n)}({\bf x}_{\A}).$$
Then we have $a_{\lambda}=0$ for all $\lambda+(d^n)\in\cP_n$ from
the linear independence of
$\{\,S^{\A}_{\lambda}\,|\,\lambda\in\cP_{\A}\,\}$ (even if it is not
a finite sum). Since $d$ is an arbitrary positive integer, it
follows that $a_{\lambda}=0$ for all $\lambda\in\cP_{\A/\B,n}$. \qed
\vskip 5mm

\section{Insertion scheme}
In this section, we will describe the combinatorial behavior of
$\A/\B$-semistandard tableaux, which is closely related to that of
rational $[n]$-semistandard tableaux. We will introduce an algorithm
of inserting an $\A/\B$-semistandard tableau into another, and
derive analogues of Robinson-Schensted-Knuth correspondence and
Littlewood-Richardson rule. We also obtain a Jacobi-Trudi type
character formula for $SST_{\A/\B}(\lambda)$.

\subsection{Robinson-Schensted-Knuth correspondence}
Let
\begin{equation}
\mathcal{F}_{\A/\B}=\bigsqcup_{c\in\mathbb{Z}}SST_{\A/\B}(c)
\end{equation}
be the set of all $\A/\B$-semistandard tableaux of level 1.

Note that for $(w^+,w^-)\in \mathcal{F}_{\A/\B}$, $w^+$ (resp.
$w^-$) is a semistandard tableau of a single row, and ${\rm
sh}(w^{\pm})\in\mathbb{Z}_{\geq 0}$, and $(w^+,w^-)\in
SST_{\A/\B}(c)$ if and only if ${\rm sh}(w^+)-{\rm
sh}(w^-)=c$.\vskip 5mm

\begin{thm}\label{RSK} For $n\geq 1$, there exists a bijection
\begin{equation*}
\kappa_{\A/\B} : \mathcal{F}_{\A/\B}^n \longrightarrow
\bigsqcup_{\lambda\in\cP_{\A/\B,n}}SST_{\A/\B}(\lambda)\times
SST_{[n]}(\lambda),
\end{equation*}
where $\mathcal{F}_{\A/\B}^n$ is the set of all $n$-tuples of
$\A/\B$-semistandard tableaux of level 1.
\end{thm}
\pf Fix $n\geq 1$. To each ordered $n$-tuple ${\bf
w}=((w^+_i,w^-_i))_{1\leq i\leq n}\in \mathcal{F}_{\A/\B}^n$, we
will associate a pair $\kappa_{\A/\B}({\bf w})=(P_{\bf w},Q_{\bf
w})\in SST_{\A/\B}(\lambda)\times SST_{[n]}(\lambda)$ for some
$\lambda\in\cP_{\A/\B,n}$, as follows. \vskip 5mm

\noindent\textsc{Step 1.} First, we let
\begin{equation*}
T^-=[(((w_1^-)^{\pi}\leftarrow (w_2^-)^{\pi})\cdots )\leftarrow
(w_n^-)^{\pi}]^{\pi}.
\end{equation*}
Let $Q$ be the recording tableau for $((((w_1^-)^{\pi}\leftarrow
(w_2^-)^{\pi})\cdots )\leftarrow (w_n^-)^{\pi})$ (see Proposition
\ref{row insertions}). That is, $\varrho_{\rm
col}((w_1^-)^{\pi},\ldots,(w_n^-)^{\pi})=((T^-)^{\pi},Q)$. We assume
that ${\rm sh}(T^-)=(d^n)/\mu$ for some $d\geq 0$ and $\mu\in
\cP_n$. Then ${\rm sh}(Q)=\delta^n_d(\mu)$.

\vskip 3mm

\noindent\textsc{Step 2.} Next, we will define $T^+$ using $w_i^+$
($1\leq i\leq n$). Since the shape of $T^+$ must be of the form
$(\lambda+(d^n))/\mu$ for some $\lambda\in\mathbb{Z}_+^n$, we will
consider the row insertions of $S_i\ast w_i^+$'s instead of
$w_i^+$'s, where $S_i$ is an $\mathbb{N}$-semistandard tableau of a
single row such that the recording tableau of the row insertions of
$S_i$'s is a rectangular complement of $Q$. Let us explain it more
precisely.

Set
$$Q^{\vee}=\delta^n_d(Q).$$
Note that ${\rm sh}(Q^{\vee})=\mu$.

Suppose that ${\rm wt}_{[n]}(Q^{\vee})=\sum_{i=1}^n\nu_i\epsilon_i$
(or the content of $Q^{\vee}$ is $\nu=(\nu_1,\ldots,\nu_n)$).
Applying Proposition \ref{row insertions} when $\A=[n]$, there
exists a unique $(S_1,\ldots,S_n)$ such that $S_i\in
SST_{[n]}(\nu_i)$ for $1\leq i\leq n$, and
$$\varrho_{\rm row}(S_1,\ldots,S_n)=(H^{\mu},Q^{\vee})\in SST_{[n]}(\mu)\times
SST_{[n]}(\mu)_{\nu}.$$ We assume that $S_i$ is empty if $\nu_i=0$.

For $1\leq i\leq n$, we put $$U_i=S_i\ast w^+_i,$$  which is an
$[n]\ast \A$-semistandard tableau of a single row with length
$\nu_i+{\rm sh}(w^+_i)$. Applying Proposition \ref{row insertions}
once again to $U_i$'s, we have
\begin{equation*}
\varrho_{\rm row}(U_1,\ldots,U_n)=(U,U_R),
\end{equation*}
where $U\in SST_{[n]\ast\A}(\lambda+(d^n))$ and $U_R\in
SST_{[n]}(\lambda+(d^n))$ for some $\lambda\in\mathbb{Z}_+^n$.
Since $i<a$ for all $i\in [n]$ and $a\in\A$ with respect to the
linear ordering on $[n]\ast\A$, we have
\begin{equation*}
U=H^{\mu}\ast T^+,
\end{equation*}
where $T^+\in SST_{\A}((\lambda+(d^n))/\mu)$.\vskip 3mm

Now, we define
\begin{equation}
P_{\bf w}=(T^+,T^-),\ \ \ Q_{\bf w}=\sigma^{-d}(U_R).
\end{equation}
Then, we have $P_{\bf w}\in SST_{\A/\B}(\lambda)$ and $Q_{\bf w}\in
SST_{[n]}(\lambda)$. Since the correspondence ${\bf w}\mapsto
(P_{\bf w},Q_{\bf w})$ is reversible by construction,
$\kappa_{\A/\B}$ is a bijection. \qed\vskip 3mm

\begin{ex}{\rm Suppose that $\A=\{a_1<a_2<a_3<\cdots\}$ and $\B=\{b_1<b_2<b_3<\cdots\}$, where all the elements are of degree $0$.
Let ${\bf w}=((w^+_i,w^-_i))_{i=1,2}$ be given by
\begin{equation*}
\begin{split}
(w^+_1,w^-_1)& = (a_1 a_1 a_2 a_4  a_5, b_3 b_3 b_4 b_6)\in SST_{\A/\B}(1),   \\
(w^+_2,w^-_2)& = (a_1 a_3 a_6 , b_2 b_3 b_6)\in SST_{\A/\B}(0).
\end{split}
\end{equation*}
Then, $$T^-=((w^-_1)^{\pi}\leftarrow (w^-_2)^{\pi} )^{\pi}= (b_6 b_4
b_3 b_3 \leftarrow b_6 b_3 b_2)^{\pi}=
\begin{array}{ccccc}
   &  &  & b_2 & b_3  \\
  b_3 & b_3 & b_4 & b_6 & b_6
 \end{array}. $$
Since the recording tableau for $(b_6 b_4 b_3 b_3 \leftarrow b_6 b_3
b_2)$ is $Q=
\begin{array}{ccccc}
1  & 1 & 1 & 1 & 2  \\
2 & 2 &  &  &
 \end{array}$, we have
$Q^{\vee}=\delta^2_5(Q)=122$, where we put $d=5$. If we put $S_1=1$
and $S_2=11$, then $\varrho_{\rm
row}(S_1,S_2)=(H^{(3,0)},Q^{\vee})$. Put
$$U_1=1 \ast a_1 a_1 a_2 a_4 a_5, \ \ \ \ U_2=11 \ast a_1 a_3 a_6.$$
Then
\begin{equation*}
\begin{split}
(U_2\rightarrow U_1)&=H^{(3,0)}\ast T^+=
\begin{array}{cccccccc}
1  & 1 & 1 & a_1 & a_3 & a_5 & a_6   \\
 a_1 & a_1 & a_2  & a_4 &
 \end{array}, \\
(U_2\rightarrow U_1)_R &=
\begin{array}{ccccccc}
1  & 1 & 1 & 1 & 1 & 1 &  2   \\
2 &  2 & 2 & 2
 \end{array}.
\end{split}
\end{equation*}
Therefore, we have
\begin{equation*}
\begin{split}
P_{\bf w}&=\left(
\begin{array}{ccccc|ccc}
     &      &     &    a_1  & a_3 & a_5 & a_6   \\
 a_1 &  a_1 & a_2 & a_4  &
 \end{array},
\begin{array}{ccccc|}
   &  &  & b_2 & b_3  \\
  b_3 & b_3 & b_4 & b_6 & b_6
 \end{array}\ \
 \right), \\
Q_{\bf w}&=
\begin{array}{c|cc}
  & 1 & 2   \\
  -2 &
 \end{array},
\end{split}
\end{equation*}
where $(P_{\bf w},Q_{\bf w})\in SST_{\A/\B}((2,-1))\times SST_{[2]}((2,-1))$.}
\end{ex}\vskip 3mm

Let us consider the character identity associated to RSK
correspondence in Theorem \ref{RSK}. Given  ${\bf
w}=((w^+_i,w^-_i))_{1\leq i\leq n}\in \mathcal{F}_{\A/\B}^n$, we set
\begin{equation*}
\begin{split}
{\rm wt}_{\A/\B}({\bf w})&=\sum_{1\leq i\leq n}{\rm
wt}_{\A/\B}(w^+_i,w^-_i)\in P_{\A}\oplus P_{\B}, \\
{\rm wt}_{[n]}({\bf w})&= \sum_{1\leq i\leq n}m_i\epsilon_i\in
P_{[n]},
\end{split}
\end{equation*}
where $m_i={\rm sh}(w^+_i)-{\rm sh}(w^-_i)$, the shape of
$(w^+_i,w^-_i)$ for $1\leq i\leq n$. Then the character of
$\mathcal{F}^n_{\A/\B}$ is given by
\begin{equation*}
\prod_{i\in [n]}\frac{\prod_{a\in \A_1}(1+x_ax_i)\prod_{b\in
\B_1}(1+x_b^{-1}x_i^{-1})}{\prod_{a\in \A_0}(1-x_ax_i)\prod_{b\in
\B_0}(1-x_b^{-1}x_i^{-1})}.
\end{equation*}

\begin{cor}\label{character for RSK} The map $\kappa_{\A/\B} : {\bf w}\mapsto (P_{\bf w},Q_{\bf w})$ in Theorem \ref{RSK} is a bijection
preserving weights, that is,
$${\rm wt}_{\A/\B}({\bf w})={\rm wt}_{\A/\B}(P_{\bf w}), \ \
{\rm wt}_{[n]}({\bf w})={\rm wt}_{[n]}(Q_{\bf w}).$$ Hence, we
obtain the following identity;
$$
\prod_{i\in [n]}\frac{\prod_{a\in \A_1}(1+x_ax_i)\prod_{b\in
\B_1}(1+x_b^{-1}x_i^{-1})}{\prod_{a\in \A_0}(1-x_ax_i)\prod_{b\in
\B_0}(1-x_b^{-1}x_i^{-1})}=
\sum_{\lambda\in\cP_{\A/\B,n}}S^{\A/\B}_{\lambda}s_{\lambda}({\bf
x}_{[n]}).
$$
\end{cor}
\pf It follows directly from the construction (see {\sc Step 1} and
{\sc Step 2} in Theorem \ref{RSK}) that ${\rm wt}_{\A/\B}({\bf
w})={\rm wt}_{\A/\B}(P_{\bf w})$. So, it suffices to show that
\begin{equation*}
{\rm wt}_{[n]}({\bf w})={\rm wt}_{[n]}(Q_{\bf w}).
\end{equation*}
Following the notations in Theorem \ref{RSK}, we have
\begin{equation*}
\begin{split}
{\rm wt}_{[n]}(Q)&=\sum_{1\leq i\leq n}{\rm
sh}(w^-_i)\epsilon_i, \\
{\rm wt}_{[n]}(Q^{\vee})&=-{\rm
wt}_{[n]}(Q)+d(\epsilon_1+\cdots+\epsilon_n),\\
{\rm wt}_{[n]}(U_R)&={\rm wt}_{[n]}(Q^{\vee})+\sum_{1\leq i\leq
n}{\rm
sh}(w^+_i)\epsilon_i,\\
{\rm wt}_{[n]}(Q_{\bf w})&={\rm
wt}_{[n]}(U_R)-d(\epsilon_1+\cdots+\epsilon_n),
\end{split}
\end{equation*}
(recall that ${\rm wt}_{[n]}(\sigma(T))={\rm
wt}_{[n]}(T)+(\epsilon_1+\cdots+\epsilon_n)$ for $T\in
SST_{[n]}(\lambda)$). Since ${\rm wt}_{[n]}({\bf w})=\sum_{1\leq
i\leq n}({\rm sh}(w^+_i)-{\rm sh}(w^-_i))\epsilon_i$, we have ${\rm
wt}_{[n]}({\bf w})={\rm wt}_{[n]}(Q_{\bf w})$. \qed

\begin{cor}\label{restriction of RSK} For $\nu=(\nu_1,\ldots,\nu_n)\in\mathbb{Z}^n$ such that $\nu_i\in\cP_{\A/\B}$, the map
$\kappa_{\A/\B}$ also gives the following bijection
\begin{equation*}
\kappa_{\A/\B} : SST_{\A/\B}(\nu_1)\times\cdots\times
SST_{\A/\B}(\nu_n) \longrightarrow
\bigsqcup_{\mu\in\cP_{\A/\B,n}}SST_{\A/\B}(\mu)\times
SST_{[n]}(\mu)_{\nu},
\end{equation*}
where $SST_{[n]}(\mu)_{\nu}$ is the set of all rational
$[n]$-semistandard tableaux of shape $\mu$ with content $\nu$, or
weight $\sum_{i=1}^n \nu_i\epsilon_i$.
\end{cor}
\pf It follows directly from Theorem \ref{RSK} and Corollary
\ref{character for RSK}. \qed\vskip 3mm

Now, we have a Jacobi-Trudi type character formula for
$S^{\A/\B}_{\lambda}$.
\begin{thm}\label{JacobiTrudi}
Suppose that $\A$ and $\B$ are infinite sets. For
$\lambda\in\mathbb{Z}_+^n$, we have
\begin{equation*}
S^{\A/\B}_{\lambda}={\rm det}\left( S^{\A/\B}_{\lambda_i-i+j}
\right)_{1\leq i,j\leq n}.
\end{equation*}
\end{thm}
\pf Since $\A$ and $\B$ are infinite sets, we have
$\cP_{\A/\B,k}=\mathbb{Z}_+^k$ for all $k\geq 1$. Fix $n\geq 1$. For
$\lambda,\mu\in\mathbb{Z}_+^n$, we define $\lambda
> \mu$ if and only if there exists an $i$ such that
$\lambda_k=\mu_k$ for $1\leq k< i$ and $\lambda_i>\mu_i$. Then $>$
is a linear ordering on $\mathbb{Z}_+^n$, called the {\it reverse
lexicographic ordering}.

Given $\mu\in\mathbb{Z}_+^n$, put $H^{\A/\B}_{\mu}=\prod_{1\leq
i\leq n}S^{\A/\B}_{\mu_i}\in\mathbb{Z}[[{\bf x}_{\A},{\bf
x}_{\B}^{-1}]]$. By Corollary \ref{restriction of RSK}, we have
$H^{\A/\B}_{\mu}=\sum_{\substack{\lambda\in\mathbb{Z}_+^n
\\ \langle \lambda\rangle =\langle \mu \rangle}}K_{\lambda\,\mu}S^{\A/\B}_{\lambda}$,
where $K_{\lambda\,\mu}=\left| SST_{[n]}(\lambda)_{\mu} \right|$. By
Lemma \ref{shift of SST}, we have for all $d>0$,
\begin{equation}\label{Kostka}
K_{\lambda\,\mu}=K_{\lambda+(d^n)\,\mu+(d^n)},
\end{equation}
which is equal to the ordinary {\it Kostka number} of shape
$\lambda+(d^n)$ with content $\mu+(d^n)$ whenever $\lambda+(d^n)$
and $\mu+(d^n)$ are ordinary partitions. This implies that
$K_{\lambda\,\mu}$ is zero unless $\lambda\geq \mu$ (cf.\cite{Mac}).
Hence, we may write
\begin{equation}\label{HS for AB}
H^{\A/\B}_{\mu}=\sum_{\substack{\lambda\geq \mu
\\ \lambda\in\cP_n}}K_{\lambda\,\mu}S^{\A/\B}_{\lambda}
+\sum_{\substack{\lambda\geq \mu
\\ \lambda\in\mathbb{Z}_+^n\setminus
\cP_n}}K_{\lambda\,\mu}S^{\A/\B}_{\lambda}.
\end{equation}

For $\lambda,\mu\in\cP_n$, let $h_{\mu}({\bf x}_{[n]})$ (resp.
$s_{\lambda}({\bf x}_{[n]})$) be the complete symmetric polynomial
(resp. Schur polynomial) in $n$ variables corresponding to $\mu$
(resp. $\lambda$). Recall that
\begin{equation}\label{HS for n}
\begin{split}
h_{\mu}({\bf x}_{[n]})&=\sum_{\lambda\geq
\mu}K_{\lambda\,\mu}s_{\lambda}({\bf x}_{[n]}),\\
s_{\lambda}({\bf x}_{[n]})&={\rm det}\left( h_{\lambda_i-i+j}({\bf
x}_{[n]}) \right)_{1\leq i,j\leq n},
\end{split}
\end{equation}
where we assume that $h_k({\bf x}_{[n]})=0$ for $k<0$ (see
\cite{Mac}).\vskip 3mm

\noindent\textsc{Case 1.} Suppose that $\lambda\in\cP_n$ is given,
where $\lambda_i-i+j>0$ for all $1\leq i,j\leq n$.  Comparing
\eqref{HS for AB} and \eqref{HS for n}, we have
\begin{equation}\label{diff}
S^{\A/\B}_{\lambda}-{\rm det}\left( H^{\A/\B}_{\lambda_i-i+j}
\right)_{1\leq i,j\leq n}= \sum_{\substack{\nu\geq \lambda
\\ \nu\in\mathbb{Z}_+^n\setminus \cP_n}}a_{\nu}S^{\A/\B}_{\nu}
\end{equation}
for some $a_{\nu}\in\mathbb{Z}$. If we apply the same argument to
$\lambda+(d^n)$ for $d\geq 0$, we have
\begin{equation}\label{diff1}
S^{\A/\B}_{\lambda+(d^n)}-{\rm det}\left(
H^{\A/\B}_{\lambda_i-i+j+d} \right)_{1\leq i,j\leq n}=
\sum_{\substack{\nu\geq \lambda
\\ \nu\in\mathbb{Z}_+^n\setminus \cP_n}}b_{\nu}S^{\A/\B}_{\nu}
\end{equation}
for some $b_{\nu}\in\mathbb{Z}$. On the other hand, by
\eqref{Kostka}, the equation \eqref{HS for AB} still holds when we
replace $\lambda$ and $\mu$  by $\lambda+(d^n)$ and $\mu+(d^n)$
($d\geq 0$), respectively. Hence, from \eqref{diff}, we also obtain
\begin{equation}\label{diff2}
S^{\A/\B}_{\lambda+(d^n)}-{\rm det}\left(
H^{\A/\B}_{\lambda_i-i+j+d} \right)_{1\leq i,j\leq n}=
\sum_{\substack{\nu\geq \lambda
\\ \nu\in\mathbb{Z}_+^n\setminus
\cP_n}}a_{\nu}S^{\A/\B}_{\nu+(d^n)}.
\end{equation}
Comparing \eqref{diff1} and \eqref{diff2}, it follows from Corollary
\ref{linear indep of SAB'} that
\begin{itemize}
\item[(1)] $a_{\nu}=b_{\nu+(d^n)}$,

\item[(2)] $a_{\nu}=0$ whenever
$\nu+(d^n)\in\cP_n$,
\end{itemize}
where we assume that $b_{\nu}=0$ for $\nu\in\cP_n$. Since $d>0$ is
arbitrary, we conclude that $a_{\nu}=0$ for all $\nu$ such that
$\nu\geq \lambda$ and $\nu\in\mathbb{Z}_+^n\setminus \cP_n$.\vskip
3mm

\noindent\textsc{Case 2.} Suppose that $\lambda\in\mathbb{Z}_+^n$ is
given and
\begin{equation*}
S^{\A/\B}_{\lambda}-{\rm det}\left( H^{\A/\B}_{\lambda_i-i+j}
\right)_{1\leq i,j\leq
n}=\sum_{\nu\in\mathbb{Z}_+^n}c_{\nu}S^{\A/\B}_{\nu},
\end{equation*}
for some $c_{\nu}\in\mathbb{Z}$. By \eqref{Kostka}, we also have
\begin{equation*}
S^{\A/\B}_{\lambda+(d^n)}-{\rm det}\left(
H^{\A/\B}_{\lambda_i-i+j+d} \right)_{1\leq i,j\leq
n}=\sum_{\nu\in\mathbb{Z}_+^n}c_{\nu}S^{\A/\B}_{\nu+(d^n)},
\end{equation*}
for all $d\geq 0$. By Case 1, the above equation is zero if $d$ is
sufficiently large, and hence $c_{\nu}=0$ for all
$\nu\in\mathbb{Z}_+^n$. This completes the proof. \qed

\begin{ex}{\rm
Suppose that  $\A$ and $\B$ are infinite sets. By Corollary
\ref{restriction of RSK}, we have
\begin{equation*}
\begin{split}
H^{\A/B}_{(1,1)}&=S^{\A/\B}_{(1,1)}+S^{\A/\B}_{(2,0)}+S^{\A/\B}_{(3,-1)}+S^{\A/\B}_{(4,-2)}+\cdots,
\\
H^{\A/B}_{(2,0)}&=\ \ \ \ \ \ \ \
+S^{\A/\B}_{(2,0)}+S^{\A/\B}_{(3,-1)}+S^{\A/\B}_{(4,-2)}+\cdots.
\end{split}
\end{equation*}
Hence,
\begin{equation*}
\begin{split}
S^{\A/\B}_{(1,1)}&=H^{\A/B}_{(1,1)}-H^{\A/B}_{(2,0)}=S^{\A/B}_{1}S^{\A/B}_{1}-S^{\A/B}_{2}S^{\A/B}_{0}
\\ &={\rm det}\left(
            \begin{array}{cc}
              S^{\A/B}_{1} & S^{\A/B}_{2} \\
              S^{\A/B}_{0} & S^{\A/B}_{1} \\
            \end{array}
          \right).
\end{split}
\end{equation*}
 }
\end{ex}

\vskip 5mm

\subsection{Littlewood-Richardson rule} Now, we are in a position to describe the
LR rule for $\A/\B$-semistandard tableaux. For $\mu\in\cP_{\A/\B,m}$
and $\nu\in\cP_{\A/\B,n}$, we will introduce an algorithm of
inserting a bitableau ${\bf T}_1=(T_1^+,T_1^-)\in SST_{\A/\B}(\mu)$
into another bitableau ${\bf T}_2=(T_2^+,T_2^-)\in SST_{\A/\B}(\nu)$
to create a new bitableau ${\bf T}=(T^+,T^-)\in
SST_{\A/\B}(\lambda)$ for some $\lambda\in\cP_{\A/\B,m+n}$ together
with a recording tableau in ${\bf LR}^{\lambda}_{\mu\,\nu}$.

Let us give a brief sketch of our algorithm, which is very similar
to RSK correspondence in Theorem \ref{RSK}. First, we define a
$\B$-semistandard tableau $T^-$ by applying the column insertion
algorithm to $(T^-_i)^{\pi}$ ($i=1,2$). Next, we consider
$\mathbb{N}$-semistandard tableaux $U_i$ ($i=1,2$) whose recording
tableau with respect to row insertion forms a rectangular complement
to that of $(T^-_i)^{\pi}$ ($i=1,2$) with respect to column
insertion. Finally, we apply the row insertion algorithm to $U_i\ast
T^+_i$ instead of $T^+_i$ ($i=1,2$) to obtain an $\A$-semistandard
tableau $T^+$ of a skew shape. Then the pair ${\bf T}=(T^+,T^-)$
becomes an $\A/\B$-semistandard tableau of shape
$\lambda\in\mathbb{Z}_+^{m+n}$, and the recording tableau is given
as an element in ${\bf LR}^{\lambda}_{\mu\,\nu}$ corresponding to
the row insertion of $U_2\ast T^+_2$ into $U_1\ast T^+_1$.

\begin{thm}\label{LR rule for SAB}
For $\mu\in\cP_{\A/\B,m}$ and $\nu\in\cP_{\A/\B,n}$, there exists a
bijection
$$\rho_{\A/\B} : SST_{\A/\B}(\mu)\times SST_{\A/\B}(\nu) \longrightarrow
\bigsqcup_{\lambda\in\cP_{\A/\B,m+n}} SST_{\A/\B}(\lambda)\times
{\bf LR}^{\lambda}_{\mu\,\nu}.$$ In terms of characters, we have
$$S^{\A/\B}_{\mu}S^{\A/\B}_{\nu}=\sum_{\lambda\in\cP_{\A/\B,m+n}}\hat{c}^{\lambda}_{\mu\,\nu}
S^{\A/\B}_{\lambda}.$$\qed
\end{thm}
\pf Let ${\bf T}_1=(T_1^+,T_1^-)\in SST_{\A/\B}(\mu)$ and ${\bf
T}_2=(T_2^+,T_2^-)\in SST_{\A/\B}(\nu)$. We will define a pair
$$\rho_{\A/\B}({\bf T}_1,{\bf T}_2)=({\bf T}, {\bf T}_R),$$
where ${\bf T}=(T^+,T^-)\in SST_{\A/\B}(\lambda)$ and ${\bf T}_R\in
{\bf LR}^{\lambda}_{\mu\,\nu}$ for some
$\lambda\in\cP_{\A/\B,m+n}$.\vskip 3mm

\noindent{\sc Step 1}. First, let us define $T^-$. Suppose that
\begin{equation*}
{\rm sh}(T^+_1)=(\mu+(d^m))/\zeta^{(1)}, \ \ {\rm
sh}(T^+_2)=(\nu+(d^n))/\zeta^{(2)},
\end{equation*}
for some $d>0$ and $\zeta^{(1)}\in\cP_m$, $\zeta^{(2)}\in\cP_n$. We
may assume that $d$ is sufficiently large. Then
\begin{equation*}
{\rm sh}(T^-_1)=(d^m)/\zeta^{(1)}, \ \ {\rm
sh}(T^-_2)=(d^n)/\zeta^{(2)}.
\end{equation*}
We define
\begin{equation*}
T^-=((T^-_1)^{\pi}\leftarrow (T^-_2)^{\pi})^{\pi}.
\end{equation*}
Then we have ${\rm sh}(T^-)=(d^{m+n})/\zeta$ for some $\zeta\in
\cP_{m+n}$.\vskip 3mm

\noindent{\sc Step 2}. Next, let us define $T^+$. Consider
$(T^-_i)^{\sharp}$ for $i=1,2$. Set $\eta^{(i)}={\rm
sh}(T^-_i)^{\sharp}$ for $i=1,2$. Note that
\begin{equation*}
\eta^{(1)}=\delta^m_d(\zeta^{(1)})',  \ \
\eta^{(2)}=\delta^n_d(\zeta^{(2)})'.
\end{equation*}
By Theorem \ref{LR rule}, there exists a unique $(S_1,S_2)\in
SST_{[d]}(\eta^{(1)})\times SST_{[d]}(\eta^{(2)})$ such that
\begin{equation*}
(S_2\rightarrow S_1)=H^{\eta},\ \ \ \ (S_2\rightarrow S_1)_R=
((T^-_2)^{\sharp} \rightarrow (T^-_1)^{\sharp})_R,
\end{equation*}
where $\eta={\rm sh}((T^-_2)^{\sharp} \rightarrow
(T^-_1)^{\sharp})$. Since
\begin{equation*}
((T^-_2)^{\sharp} \rightarrow (T^-_1)^{\sharp})= ((T^-_1)^{\pi}
\leftarrow (T^-_2)^{\pi})^t=[(T^-)^{\pi}]^t=(T^-)^{\sharp},
\end{equation*}
we have
$$\eta=\delta^{m+n}_d(\zeta)'.$$
\vskip 3mm

Set $U_i=\delta^d_{m}(S_i)^t$ ($i=1,2$). Then
\begin{equation*}
{\rm sh}(U_i)={\rm sh}(\delta^d_{m}(S_i))'
=\delta^d_m(\delta^m_d(\zeta^{(i)})')' =\delta^d_m
(\delta^d_m((\zeta^{(i)})') )'=\zeta^{(i)},
\end{equation*}
and hence $U_i\in SST_{[d]'}(\zeta^{(i)})$ ($i=1,2$). Then
\begin{equation*}
\begin{split}
(U_2\rightarrow U_1)
&=(\delta^d_{m}(S_2)^t\rightarrow\delta^d_m(S_1)^t)
=[\delta^d_{m}(S_1)\leftarrow \delta^d_m(S_1)]^t \\
&=[\delta^d_{m+n}(S_2\rightarrow S_1)]^t \ \ \ \ \text{by Theorem
\ref{Stroomer}}
\\
&=[\delta^d_{m+n}(H^{\eta})]^t,
\end{split}
\end{equation*}
and
\begin{equation*}
\begin{split}
{\rm sh}(U_2\rightarrow U_1)
&=\delta^d_{m+n}(\eta)'=\delta^d_{m+n}(\delta_d^{m+n}(\zeta)')' \\
&= \delta^d_{m+n}(\delta^d_{m+n}(\zeta'))'=\zeta.
\end{split}
\end{equation*}
Hence, $(U_2\rightarrow U_1)\in SST_{[d]'}(\zeta)$.

Now, we set
\begin{equation*}
\begin{split}
& \widehat{U}_1=U_1\ast T^+_1\in SST_{[d]'\ast \A}(\mu+(d^m)), \\
& \widehat{U}_2=U_2\ast T^+_2\in SST_{[d]'\ast \A}(\nu+(d^n)).
\end{split}
\end{equation*}
Then, we have
\begin{equation*}
(\widehat{U}_2\rightarrow \widehat{U}_1)=(U_2\rightarrow U_1)\ast
T^+,
\end{equation*}
where $T^+\in SST_{\A}((\lambda+(d^{m+n}))/\zeta)$ for some
$\lambda\in\mathbb{Z}_+^{m+n}$. We define
\begin{equation}
{\bf T}=(T^+,T^-).
\end{equation}
We can check that ${\bf T}$ does not depend on the choice of $d$,
and ${\bf T}\in SST_{\A/\B}(\lambda)$. \vskip 3mm

\noindent{\sc Step 3.} By Lemma \ref{LR tableaux}, we have
$((\widehat{U}_2\rightarrow \widehat{U}_1)_R)^t\in
LR^{(\lambda+(d^{m+n}))'}_{(\mu+(d^m))'\,(\nu+(d^n))'}$. Now, we
define
\begin{equation}
{\bf T}_R=[Q]\in {\bf LR}^{\lambda}_{\mu\,\nu},
\end{equation}
where $[Q]$ is the element in ${\bf LR}^{\lambda}_{\mu\,\nu}$
including $((\widehat{U}_2\rightarrow \widehat{U}_1)_R)^t$ (see
\eqref{skew LR tableaux for rational}).\vskip 3mm

Since the construction of $({\bf T},{\bf T}_R)$ is reversible,
$\rho_{\A/\B}$ is a bijection. This completes the proof. \qed

\begin{ex}{\rm
Suppose that $\A=\{a_1<a_2<a_3<\cdots\}$ and
$\B=\{b_1<b_2<b_3<\cdots\}$. For convenience, we assume that all the
elements are of degree $0$. Suppose that
\begin{equation*}
\begin{split}
{\bf T}_1&=(T_1^+,T_1^-)= \left(
\begin{array}{c|cc}
  a_1 & a_2 & a_2 \\
   &  &
\end{array}
,
\begin{array}{c|cc}
  b_2 &  &  \\
  b_3 &  &
\end{array}
\right)\in SST_{\A/\B}((2,-1)),
 \\
{\bf T}_2&=(T_2^+,T_2^-)=\left(
\begin{array}{cc|c}
   & a_3 & a_4 \\
  a_1 &  &
\end{array}
,
\begin{array}{cc|c}
   & b_1 &  \\
  b_1 & b_2 &
\end{array}
\right)\in SST_{\A/\B}((1,-1)).
\end{split}
\end{equation*}
Then we have
$$T^-=
\begin{array}{cc|}
   &  \\
   & b_1 \\
  b_1 & b_2 \\
  b_2 & b_3
\end{array}\ \ .
$$
Note that
$$
(T^-_1)^{\sharp}=
\begin{array}{cc}
  b_3 & b_2   \\
   &
\end{array},\ \ \
(T^-_2)^{\sharp}=
\begin{array}{cc}
  b_2 & b_1   \\
  b_1 &
\end{array},
$$
and
$$
(T^-_2)^{\sharp}\rightarrow (T^-_1)^{\sharp}=
\begin{array}{cccc}
  b_3 & b_2 & b_1 & \\
  b_2 & b_1 & &
\end{array},  \ \
((T^-_2)^{\sharp}\rightarrow (T^-_1)^{\sharp})_R=
\begin{array}{cccc}
  \bullet & \bullet & 1 & \\
  1 & 2 & &
\end{array}.
$$
Hence, if we put
$S_1=\begin{array}{cc}
       1 & 2
     \end{array}
$ and $S_2=\begin{array}{cc}
             1 & 1 \\
             2 &
           \end{array}
$, then $S_2\rightarrow S_1 = H^{(3,2)}$ and $(S_2\rightarrow
S_1)_R=((T^-_2)^{\sharp}\rightarrow (T^-_1)^{\sharp})_R$. Now, we
put
$$U_1=\delta^4_2(S_1)^t=
\begin{array}{ccc}
  1 & 3 & 4 \\
  2 & 3 & 4
\end{array}, \ \ U_2=\delta^4_2(S_2)^t=
\begin{array}{ccc}
  2 & 3 & 4 \\
  3 & 4 &
\end{array}.$$
and
\begin{equation*}
\begin{split}
\widehat{U}_1&=U_1\ast T^+_1=
\begin{array}{cccccc}
  1 & 3 & 4 & a_1 & a_2 & a_2\\
  2 & 3 & 4
\end{array}, \\
\widehat{U}_2&=U_2\ast T^+_2=
\begin{array}{ccccc}
  2 & 3 & 4 & a_3 & a_4\\
  3 & 4 &   a_1
\end{array}.
\end{split}
\end{equation*}
Then we have
\begin{equation*}
\begin{split}
\widehat{U}_2 \rightarrow \widehat{U}_1&=
\begin{array}{cccccccc}
  1 & 2 & 3 & 4 & a_1 & a_2 & a_3 & a_4 \\
  2 & 3 & 4 & a_1 &   &   &   &   \\
  3 & 4 & a_2 &   &   &   &   &   \\
  3 & 4 &   &   &   &   &   &
\end{array},\\
(\widehat{U}_2 \rightarrow \widehat{U}_1)_R&=
\begin{array}{cccccccc}
  \bullet  &  \bullet &  \bullet &  \bullet &  \bullet &  \bullet & 4 & 5 \\
  \bullet  &  \bullet &  \bullet & 3 &   &   &   &   \\
  1 & 2 & 3 &   &   &   &   &   \\
  1 & 2 &   &   &   &   &   &
\end{array}.
\end{split}
\end{equation*}
Note that $((\widehat{U}_2\rightarrow \widehat{U}_1)_R)^t\in
LR^{(8,4,3,2)}_{(6,3)\,(5,3)}$. Therefore,
\begin{equation*}
{\bf T}=(T^+,T^-)=\left(
\begin{array}{cccc|cccc}
    &   &   &   & a_1 & a_2 & a_3 & a_4 \\
    &   &   & a_1 &   &   &   &   \\
    &   & a_2 &   &   &   &   &   \\
    &   &   &   &   &   &   &
\end{array},
\begin{array}{cc|}
   &  \\
   & b_1 \\
  b_1 & b_2 \\
  b_2 & b_3
\end{array}\ \ \ \ \right),
\end{equation*}
where ${\bf T} \in SST_{\A/\B}(4,0,-1,-2)$, and $${\bf T}_R=\left[
\begin{array}{ccc}
  \bullet  &  \bullet &  1 \\
  \bullet  & 1 &    \\
  \bullet  &  &    \\
  \bullet  &  &    \\
  2 &  &      \\
  3 &  &
\end{array}\right]\in {\bf LR}^{(4,0,-1,-2)}_{(2,-1)(1,-1)},$$
(see Lemma \ref{LR coeff for SAB}).}
\end{ex}\vskip 5mm

\subsection{Skew Littlewood-Richardson rule}
We have seen that the products of characters of
$SST_{\A/\B}(\lambda)$'s is determined by the branching coefficients
of rational Schur polynomials. This might be understood as a dual
relation between $\A/\B$-semistandard tableaux and rational
$[n]$-semistandard tableaux. To complete this dual relationship, we
will define $\A/\B$-semistandard tableaux of skew shapes, and then
obtain their skew LR rule, which is completely determined by the LR
rule of rational $[n]$-semistandard tableaux.

\begin{df}\label{skew ABsst}{\rm
Suppose that $\lambda$ and $\mu\in\mathbb{Z}_+^n$ are given.  An
{\it $\A/\B$-semistandard tableau of skew shape $\lambda/\mu$} is a
pair of tableaux $(T^+,T^-)$ such that
$$T^+\in SST_{\A}((\lambda+(d^n))/\nu), \ \ \ \ T^- \in SST_{\B}((\mu+(d^{n}))/\nu),$$
for some integer $d\geq 0$ and $\nu\in\cP_n$ satisfying
\begin{itemize}
\item[(1)] $\lambda+(d^n), \mu+(d^n) \in\cP_n$,

\item[(2)] $\nu\subset (d^n), \lambda+(d^n)$, and $\mu+(d^n)$.
\end{itemize}}
\end{df}\vskip 3mm
We denote by $SST_{\A/\B}(\lambda/\mu)$ the set of all
$\A/\B$-semistandard tableaux of skew shape $\lambda/\mu$. Note that
$SST_{\A/\B}(\lambda/{\bf 0}_n)=SST_{\A/\B}(\lambda)$. The character
of $SST_{\A/\B}(\lambda/\mu)$ is defined similarly, and denoted by
$S^{\A/\B}_{\lambda/\mu}$.\vskip 3mm

To discuss the skew LR rule, we need to consider a rectangular
complement  of a Littlewood-Richardson tableau, which is given in
the following lemma.

\begin{lem}\label{dual LR}
For $\lambda,\mu$, and $\nu\in\cP_n$, there exists a bijection
$$\delta^n_{p,q} : LR^{\lambda}_{\mu\,\nu}\longrightarrow LR^{\delta^n_{p+q}(\lambda)}_{\delta^n_p(\mu)\,\delta^n_q(\nu)},$$
for $p\geq \mu_1$ and $q\geq \nu_1$.
\end{lem}
\pf Suppose that $LR^{\lambda}_{\mu\,\nu}$ is non-empty, and  $Q\in
LR^{\lambda}_{\mu\,\nu}$ is given. Then $\tau(Q)\in
LR^{\lambda'}_{\mu'\,\nu'}$, where $\tau$ is the map given in
Corollary \ref{conjugate LR}. By Theorem \ref{LR rule}, there exists
a unique pair $(T_1,T_2)\in SST_{[n]}(\mu)\times SST_{[n]}(\nu)$
such that $(T_2\rightarrow T_1)=H^{\lambda}$ and $((T_2\rightarrow
T_1)_R)^t=\tau(Q)$. Now, we define
\begin{equation}
\delta^n_{p,q}(Q)=(\delta^n_{p}(T_1)\leftarrow \delta^n_{q}(T_2))_R,
\end{equation}
where $p\geq \mu_1$ and $q\geq \nu_1$. By Theorem \ref{Stroomer},
$(\delta^n_{p}(T_1)\leftarrow
\delta^n_{q}(T_2))=\delta^n_{p+q}(T_2\rightarrow
T_1)=\delta^n_{p+q}(H^{\lambda})$, and hence $\delta^n_{p,q}(Q)\in
LR^{\delta^n_{p+q}(\lambda)}_{\delta^n_p(\mu)\,\delta^n_q(\nu)}$.

Suppose that $\delta^n_{p,q}(Q)=\delta^n_{p,q}(Q')$ for some $Q'\in
LR^{\lambda}_{\mu\,\nu}$. Let $(T'_1,T'_2)\in SST_{[n]}(\mu)\times
SST_{[n]}(\nu)$ be the associated pair such that $(T'_2\rightarrow
T'_1)=H^{\lambda}$ and $((T'_2\rightarrow T'_1)_R)^t=\tau(Q')$.

Note that
\begin{equation*}
\begin{split}
(\delta^n_{p}(T'_1)\leftarrow
\delta^n_{q}(T'_2))&=\delta^n_{p+q}(T'_2\rightarrow
T'_1)=\delta^n_{p+q}(H^{\lambda})=(\delta^n_{p}(T_1)\leftarrow
\delta^n_{q}(T_2)), \\
(\delta^n_{p}(T'_1)\leftarrow \delta^n_{q}(T'_2))_R&=
\delta^n_{p,q}(Q')=\delta^n_{p,q}(Q)=(\delta^n_{p}(T_1)\leftarrow
\delta^n_{q}(T_2))_R.
\end{split}
\end{equation*}
Since $\rho_{\rm col}$ is a bijection, we have
$\delta^n_{p}(T'_1)=\delta^n_{p}(T_1)$ and
$\delta^n_{q}(T'_2)=\delta^n_{q}(T_2)$. Since $\delta^n_k$ is also
bijective for $k=p,q$, we have $(T_1,T_2)=(T'_1,T'_2)$, and hence
$Q=Q'$, which implies that $\delta^n_{p,q}$ is one-to-one.

Since $\delta^n_{p,q}$ also gives a one-to-one map from
$LR^{\delta^n_{p+q}(\lambda)}_{\delta^n_p(\mu)\,\delta^n_q(\nu)}$ to
$LR^{\lambda}_{\mu\,\nu}$, $\delta^n_{p,q}$ is a bijection.
\qed\vskip 3mm

\begin{thm}\label{skew LR for SAB}
For $\lambda, \mu\in\mathbb{Z}_+^n$, there exists a bijection
$$J_{\A/\B} : SST_{\A/\B}(\lambda/\mu)\longrightarrow \bigsqcup_{\nu\in\cP_{\A/\B,n}}
SST_{\A/\B}(\nu)\times {\bf LR}^{\lambda/\mu}_{\nu}.$$ In terms of
characters, we have
$$S^{\A/\B}_{\lambda/\mu}=\sum_{\nu\in\cP_{\A/\B,n}}c^{\lambda}_{\mu\,\nu}S^{\A/\B}_{\nu}.$$
\end{thm}
\pf To each ${\bf T}=(T^+,T^-)\in SST_{\A/\B}(\lambda/\mu)$, we will
associate a pair $J_{\A/\B}({\bf T})=(j({\bf T}),j({\bf T})_R)\in
SST_{\A/\B}(\nu)\times {\bf LR}^{\lambda/\mu}_{\nu}$ for some
$\nu\in\cP_{\A/\B,n}$.\vskip 3mm

First, consider $(T^-)^{\pi}$. Suppose that ${\rm
sh}((T^-)^{\pi})=\alpha/\beta$ for some $\alpha,\beta\in \cP_n$. In
fact, one may assume that $\beta=\delta^n_{r}(\mu+(p^n))$, where
$p=\max\{-\mu_n,0\}$ and $r=\mu_1+p$.

Applying Corollary \ref{skew LR} to $(T^-)^{\pi}$, we have
\begin{equation*}
J((T^-)^{\pi})=(j((T^-)^{\pi}),j((T^-)^{\pi})_R)\in
SST_{\B^{\pi}}(\gamma)\times LR^{\alpha}_{\beta\,\gamma},
\end{equation*}
for some $\gamma\in\cP_n$. We define
\begin{equation*}
\widehat{T}^-=j((T^-)^{\pi})^{\pi}.
\end{equation*}
Note that ${\rm sh}(\widehat{T}^-)=\gamma^{\pi}$.

Next, consider $Q=j((T^-)^{\pi})_R$. If we put $q=\gamma_1$, then by
Lemma \ref{dual LR} we have
$$\delta^n_{p,q}(Q)=Q^{\vee}\in
LR^{\delta^n_{r+q}(\alpha)}_{\delta^n_r(\beta)\,\delta^n_q(\gamma)}.$$
By definition of $\alpha$ and $\beta$, we can check the following
facts;
\begin{itemize}
\item[(1)] ${\rm
sh}(T^+)=[\lambda+((p+q)^n)]/\delta^n_{r+q}(\alpha)$,

\item[(2)] $\delta^n_r(\beta)=\mu+(p^n)$,

\item[(3)] $Q^{\vee}\ast T^+\in
SST_{\mathbb{N}\ast\A}([\lambda+((p+q)^n)]/\mu+(p^n))$.
\end{itemize}
Since $j(Q^{\vee})=H^{\delta^n_q(\gamma)}$ by Theorem
\ref{switching} (3), we obtain
\begin{equation*}
J(Q^{\vee}\ast T^+)=(H^{\delta^n_q(\gamma)}\ast
\widehat{T}^+,\widehat{Q}),
\end{equation*}
where
\begin{itemize}
\item[(1)] $\widehat{T}^+\in SST_{\A}((\nu+(q^n))/\delta^n_q(\gamma))$ for some
$\nu\in\mathbb{Z}_+^n$,

\item[(2)] $\widehat{Q}\in LR^{\lambda+((p+q)^n)}_{\mu+(p^n)\,\nu+(q^n)}$.
\end{itemize}
Now, we define
\begin{equation}
j({\bf T})=(\widehat{T}^+,\widehat{T}^-).
\end{equation}
Then $j({\bf T})\in SST_{\A/\B}(\nu)$ since ${\rm
sh}(\widehat{T}^-)=\gamma^{\pi}=(q^n)/\delta^n_q(\gamma)$. And we
define $j({\bf T})_R$ to be the element in ${\bf
LR}^{\lambda/\mu}_{\nu}$ containing $\widehat{Q}$ (see \eqref{LR
tableaux for rational}). \vskip 3mm

Since our construction is reversible, the correspondence ${\bf
T}\mapsto (j({\bf T}),j({\bf T})_R)$ is bijective. Moreover, we
obtain the corresponding identity from the characters of the both
sides since ${\rm wt}_{\A/\B}({\bf T})={\rm wt}_{\A/\B}(j({\bf T}))$
for all $T\in SST_{\A/\B}(\lambda/\mu)$. This completes the proof.
\qed

\begin{cor} For $\lambda\in\mathbb{Z}_+^n$, we have
$$S^{\A/\B}_{{\bf 0}_n/\lambda}=S^{\A/\B}_{\lambda^*}.$$
\end{cor}
\pf By Theorem \ref{skew LR for SAB}, it suffices to show that
\begin{equation*}
c^{{\bf 0}_n}_{\lambda\,\mu}=
\begin{cases}
1, & \text{if $\mu=\lambda^*$}, \\
0, & \text{otherwise}.
\end{cases}
\end{equation*}
Note that $c^{{\bf
0}_n}_{\lambda\,\mu}=N^{((p+q)^n)}_{\lambda+(p^n)\, \mu+(q^n)}$ for
sufficiently large $p,q>0$ whenever $\lambda+(p^n),
\mu+(q^n)\in\cP_n$ and $\lambda+(p^n)\subset ((p+q)^n)$. Fix such
$p$ and $q$. Then it is not difficult to see that there exists a
unique Littlewood-Richardson tableau $Q$ of shape
$(p+q)^n/(\lambda+(p^n))$, whose content should be
$\delta^n_{p+q}(\lambda+(p^n))$. Hence, we have
\begin{equation*}
\mu+(q^n)=\delta^n_{p+q}(\lambda+(p^n))=\lambda^*+(q^n),
\end{equation*}
which implies that $\mu=\lambda^*$. \qed

\section{Representations of infinite dimensional Lie superalgebras}

In this section, we discuss a relation between $\A/\B$-semistandard
tableaux and a certain class of representations of infinite
dimensional Lie superalgebra $\widehat{\frak{gl}}_{\infty|\infty}$
(or $\widehat{\frak{gl}}_{\infty}$). More precisely, we will show
that the characters of certain quasi-finite irreducible
representations of $\widehat{\frak{gl}}_{\infty|\infty}$ (or
$\widehat{\frak{gl}}_{\infty}$) parameterized by generalized
partitions (see \cite{CL,KacR2}) are realized as those of
$\A/\B$-semistandard tableaux of the corresponding shapes with
suitable choices of $\A$ and $\B$. Using the combinatorial results
established in the previous sections, we will characterize the
Grothendieck rings for certain categories of semi-simple
representations of $\widehat{\frak{gl}}_{\infty|\infty}$ (or
$\widehat{\frak{gl}}_{\infty}$), whose irreducible factors are
parameterized by generalized partitions.  \vskip 3mm

\subsection{Fock space representations  }

Let $\frac{1}{2}\mathbb{Z}=\{\,\frac{n}{2}\,|\,n\in\mathbb{Z}\,\}$
be a $\mathbb{Z}_2$-graded set with
$\left(\frac{1}{2}\mathbb{Z}\right)_0=\mathbb{Z}$ and
$\left(\frac{1}{2}\mathbb{Z}\right)_1=\frac{1}{2}+\mathbb{Z}$. A
linear ordering on $\frac{1}{2}\mathbb{Z}$ is given as the ordinary
one. Let $\mathbb{C}^{\infty|\infty}$ be the associated superspace
with a basis $\{\,\epsilon_k\,|\,k\in\frac{1}{2}\mathbb{Z}\,\}$. Let
\begin{equation}
{\frak{gl}}_{\infty|\infty}=\{\,(a_{ij})_{i,j\in
\frac{1}{2}\mathbb{Z}}\,|\,a_{ij}\in\mathbb{C},\ a_{ij}=0\ \
\text{for $|i-j|\gg 0$}\,\}.
\end{equation}
Since an element in ${\frak{gl}}_{\infty|\infty}$ is a linear
transformation of $\mathbb{C}^{\infty|\infty}$,
${\frak{gl}}_{\infty|\infty}$ is naturally endowed with a
$\mathbb{Z}_2$-grading, and becomes a Lie superalgebra with respect
to super commutator. For $i,j\in \frac{1}{2}\mathbb{Z}$, we denote
by $e_{ij}$ the elementary matrix with 1 at the $i$th row and the
$j$th column and $0$ elsewhere.

Let
$\widehat{\frak{gl}}_{\infty|\infty}={\frak{gl}}_{\infty|\infty}\oplus
\mathbb{C}K$ be a central extension of ${\frak{gl}}_{\infty|\infty}$
with respect to the following two-cocycle
\begin{equation*}
\alpha(A,B)={\rm str}([J,A]B) \ \ \ \ \text{($A,B\in
{\frak{gl}}_{\infty|\infty}$)},
\end{equation*}
where $J=\sum_{r\leq 0}e_{rr}$, and ${\rm str}$ is the supertrace
defined by ${\rm str}((a_{ij}))=\sum_{i\in
\frac{1}{2}\mathbb{Z}}(-1)^{2i}a_{ii}$. Then we have a triangular
decomposition
\begin{equation*}
\widehat{\frak{gl}}_{\infty|\infty}=\frak{n}_+ \oplus \frak{h}
\oplus \frak{n}_-,
\end{equation*}
where $\frak{h}$ is the subalgebra spanned by diagonal matrices and
$K$, and $\frak{n}^+$ (resp. $\frak{n}^-$) is the subalgebra of
strictly upper (resp. lower) triangular matrices. With this, one can
define a Verma module $M(\Lambda)$ of
$\widehat{\frak{gl}}_{\infty|\infty}$ with highest weight
$\Lambda\in\frak{h}^*$. Then we denote by $L(\Lambda)$ the unique
irreducible quotient of $M(\Lambda)$ with highest weight $\Lambda$.
If we define ${\rm deg}e_{ij}=j-i$ for
$i,j\in\frac{1}{2}\mathbb{Z}$, then
$\widehat{\frak{gl}}_{\infty|\infty}$ becomes a
$\frac{1}{2}\mathbb{Z}$-graded Lie superalgebra. And if we define
the degree of the highest weight vector in $L(\Lambda)$ to be $0$,
then $L(\Lambda)$ is also naturally $\frac{1}{2}\mathbb{Z}$-graded
$L(\Lambda)=\bigoplus_{k\in \frac{1}{2}\mathbb{Z}}L(\Lambda)_k$. We
say that $L(\Lambda)$ is {\it quasi-finite} if ${\rm
dim}L(\Lambda)_k$ is finite for all $k\in \frac{1}{2}\mathbb{Z}$
(cf.\cite{CW03,KacR1}).

For $n\geq 1$, let $\mathfrak{F}^n$ be the infinite dimensional Fock
space generated by $n$ pairs of free fermions and $n$ pairs of free
bosons (see \cite{CL,CW03,KacL} for a detailed description). Then we
have a natural commuting action of
$\widehat{\frak{gl}}_{\infty|\infty}$ and $\frak{gl}_n$ on
$\mathfrak{F}^n$. Using Howe duality, Cheng and Wang proved the
following multiplicity-free decomposition of $\mathfrak{F}^n$.
\begin{thm}[\cite{CW03}]\label{CW decomposition}
As a $(\widehat{\frak{gl}}_{\infty|\infty},\frak{gl}_n)$-module,
$$\mathfrak{F}^n\simeq \bigoplus_{\lambda\in\mathbb{Z}_+^n}L(\Lambda(\lambda))\otimes L_n(\lambda),$$
where $L_n(\lambda)$ is the irreducible rational representation of
$\frak{gl}_n$ corresponding to $\lambda$, and
$\Lambda(\lambda)\in\frak{h}^*$ is the highest weight determined by
\begin{equation*}
\begin{split}
\Lambda(\lambda)(e_{kk})&=
\begin{cases}
\max\{\lambda'_k-k,0\}, & \text{if $k\in\mathbb{Z}_{> 0}$}, \\
-\max\{\lambda'_{k-1}+k,0\}, & \text{if $k\in\mathbb{Z}_{\leq 0}$}, \\
\max\{\lambda_{k+\frac{1}{2}}-k+\frac{1}{2},0\}, & \text{if $k\in\frac{1}{2}+\mathbb{Z}_{\geq 0}$}, \\
-\max\{-\lambda_{n+k+\frac{1}{2}}+k-\frac{1}{2},0\}, & \text{if $k\in -\frac{1}{2}-\mathbb{Z}_{\geq 0}$}, \\
\end{cases} \\
\Lambda(\lambda)(K)&=n,
\end{split}
\end{equation*}
{\rm (} $\lambda'_i$  is the number of nodes in the $i$th column of
$\lambda$ for $i\in\mathbb{Z}\setminus\{0\}$ {\rm )}.\qed
\end{thm}

For $k\in\frac{1}{2}\mathbb{Z}$, let $\omega_k$ be the fundamental
weight given by $\omega_k(e_{ll})=\delta_{kl}$ and $\omega_k(K)=0$,
and let $x_k=e^{\omega_k}$ be the formal variable. Then we can
define the character of $L(\Lambda(\lambda))$
$(\lambda\in\mathbb{Z}_+^n)$ with respect to the action of the
abelian subalgebra
$\bigoplus_{k\in\frac{1}{2}\mathbb{Z}}\mathbb{C}e_{kk}$. In
\cite{CL}, using the classical Cauchy identities of (hook) Schur
functions (cf.\cite{Mac,Rem}), Cheng and Lam showed that
$L(\Lambda(\lambda))$ is given as a linear combination of product of
two hook Schur functions.

Put $\A=\left(\frac{1}{2}\mathbb{Z}_{>0}\right)'$ and
$\B=\left(\frac{1}{2}\mathbb{Z}_{\leq 0}\right)'$. Then we may view
${\bf x}_{\A}=\{\,e^{\omega_k}\,|\,k\in\frac{1}{2}\mathbb{Z}_{>
0}\,\}$ and ${\bf
x}_{\B}=\{\,e^{\omega_k}\,|\,k\in\frac{1}{2}\mathbb{Z}_{\leq
0}\,\}$. Now, we obtain a new combinatorial realization of ${\rm
ch}L(\Lambda(\lambda))$.
\begin{thm}\label{CL formula}
For $\lambda\in\mathbb{Z}_+^n$, we have
$${\rm ch}L(\Lambda(\lambda))=S^{\A/\B}_{\lambda},$$
where $\A=\left(\frac{1}{2}\mathbb{Z}_{>0}\right)'$ and
$\B=\left(\frac{1}{2}\mathbb{Z}_{\leq 0}\right)'$.
\end{thm}
\pf From the Cheng and Lam's formula (Theorem 3.2 in \cite{CL}), we
have
\begin{equation*}
{\rm
ch}L(\Lambda(\lambda))=\sum_{\mu,\nu\in\cP_n}c^{\lambda}_{\mu\,\nu^*}
S_{\mu}({\bf x}_{{\A_0\ast\A_1}})S_{\nu}({\bf
x}_{{\B_0\ast\B_1}}^{-1}).
\end{equation*}
Note that $S_{\mu}({\bf x}_{{\A_0\ast\A_1}})$ and $S_{\nu}({\bf
x}_{{\B_0\ast\B_1}}^{-1})$ are hook Schur functions with countably
many even and odd variables. By Lemma \ref{switching character}, we
have $S^{\A_0\ast\A_1}_{\mu}=S^{\A}_{\mu}$ and
$S^{\B_0\ast\B_1}_{\nu}=S^{\B}_{\nu}$. Therefore, ${\rm
ch}L(\Lambda(\lambda))=S^{\A/\B}_{\lambda}$ by Proposition
\ref{character for SAB}. \qed\vskip 3mm

\begin{rem}{\rm It is not difficult to see that
the left-hand sides in the character identities of RSK
correspondence in Theorem \ref{RSK} (or Corollary \ref{character for
RSK}) and the Fock space decomposition in Theorem \ref{CW
decomposition} are equal. Comparing these two identities, we can
also prove that ${\rm ch}L(\Lambda(\lambda))=S^{\A/\B}_{\lambda}$
from the linear independence of rational Schur polynomials.}
\end{rem}

Now, we have a Jacobi-Trudi type character formula for
$L(\Lambda(\lambda))$.
\begin{cor}
Under the above hypothesis, we have
$${\rm
ch}L(\Lambda(\lambda))={\rm det}\left( {\rm
ch}L(\Lambda(\lambda_i-i+j)) \right)_{1\leq i,j\leq n}.$$
\end{cor}
\pf It follows directly from Theorem \ref{JacobiTrudi}. \qed

\begin{rem}{\rm
(1) For $\lambda\in\mathbb{Z}_+^n$, the tableau ${\bf T}^{\lambda}$
in $SST_{\A/\B}(\lambda)$ corresponding to the highest weight vector
can be found easily. For example, if $\lambda=(4,3,2,-2,-3)$, then
${\bf T}^{\lambda}$ is given by filling the generalized Young
diagram $\lambda$ in the following pattern \vskip 3mm
\begin{center}
${\bf T}^{\lambda}=$
\begin{tabular}{ccc|cccc}
   & & & $\frac{1}{2}$ & $\frac{1}{2}$ & $\frac{1}{2}$ & $\frac{1}{2}$ \\
   & & & $1$ & $\frac{3}{2}$ & $\frac{3}{2}$ &  \\
   & &  & $1$ & $2$ & & \\
   & $-1$ & $0$ & & & &  \\
   $-\frac{1}{2}$ & $-\frac{1}{2}$ & $0$ & & & &
\end{tabular}.
\end{center}\vskip 3mm
In this case, we have
$$\Lambda(\lambda)=\omega_2+2\omega_{\frac{3}{2}}+2\omega_1
+4\omega_{\frac{1}{2}}-2\omega_0-2\omega_{-\frac{1}{2}}-\omega_{-1}+5\Lambda_0,$$
where $\Lambda_0\in\frak{h}^*$ is defined by $\Lambda_0(e_{kk})=0$
for all $k\in\frac{1}{2}\mathbb{Z}$, and $\Lambda_0(K)=1$.

(2) The LR rule for $\A/\B$-semistandard tableaux (Theorem \ref{LR
rule for SAB}) gives a combinatorial interpretation of decomposition
of the tensor product $L(\Lambda(\lambda))\otimes L(\Lambda(\mu))$
(cf.Theorem 6.1 in \cite{CL}). In particular, the RSK correspondence
(Theorem \ref{RSK}) corresponds to the decomposition of $\frak{F}^n$
as a $(\widehat{\frak{gl}}_{\infty|\infty},\frak{gl}_n)$-module
(Theorem \ref{CW decomposition}).

(3) It would be interesting to construct an explicit basis of
$L(\Lambda(\lambda))$ whose elements are parameterized by
$\A/\B$-semistandard tableaux of shape $\lambda$. }
\end{rem}\vskip 5mm

Next, let ${\frak{gl}}_{\infty}$ be the subalgebra of
${\frak{gl}}_{\infty|\infty}$ consisting of matrices
$(a_{ij})_{i,j\in\frac{1}{2}\mathbb{Z}}$ such that $a_{ij}=0$ unless
$i,j\in\mathbb{Z}$. Put
$\widehat{{\frak{gl}}}_{\infty}={\frak{gl}}_{\infty}\oplus
\mathbb{C}K$, which is an infinite dimensional Lie algebra. The
triangular decomposition is naturally induced from
$\widehat{\frak{gl}}_{\infty|\infty}$, say,
${\frak{gl}}_{\infty}=\frak{n}_+^0\oplus \frak{h}^0\oplus
\frak{n}_+^0$. As in the case of
$\widehat{\frak{gl}}_{\infty|\infty}$, we can define the Verma
module $M^0(\Lambda)$, and the associated irreducible highest weight
module $L^0(\Lambda)$ for $\Lambda\in(\frak{h}^0)^*$. The
fundamental weights $\omega_k$ ($k\in\mathbb{Z}$) and $\Lambda_0$
are still available.

For $n\geq 1$, let $\mathfrak{F}_0^n$ be the infinite dimensional
Fock space generated by $n$ pairs of free bosons (see \cite{KacR2}
for a detailed description). In \cite{KacR2}, using a natural
commuting action of $\widehat{\frak{gl}}_{\infty|\infty}$ and
$\frak{gl}_n$ on $\mathfrak{F}_0^n$, Kac and Radul derived a
multiplicity-free decomposition as follows;
\begin{thm}[\cite{KacR2}]\label{KR decomposition}
As a $(\widehat{\frak{gl}}_{\infty},\frak{gl}_n)$-module,
$$\mathfrak{F}_0^n\simeq \bigoplus_{\lambda\in\mathbb{Z}_+^n}L^0(\Lambda(\lambda))\otimes L_n(\lambda),$$
where $\Lambda(\lambda)\in(\frak{h}^0)^*$ is the highest weight
determined by
\begin{equation*}
\begin{split}
\Lambda(\lambda)(e_{kk})&=
\begin{cases}
\lambda_k, & \text{if $k\in\mathbb{Z}_{> 0}$ and $\lambda_k > 0$}, \\
\lambda_{n+k}, & \text{if $k\in\mathbb{Z}_{\leq 0}$ and
$\lambda_{n+k}<0$}, \\
0, & \text{otherwise,}
\end{cases} \\
\Lambda(\lambda)(K)&=-n.
\end{split}
\end{equation*}
\qed
\end{thm}

Let $\A=\mathbb{Z}_{>0}$ and $\B=\mathbb{Z}_{\leq 0}$ be the sets
with the usual linear ordering such that all the elements are of
degree 0, that is, $\A_0=\A$ and $\B_0=\B$. Comparing the Kac and
Radul's formula for $L^0(\Lambda(\lambda))$ with Proposition
\ref{character for SAB}, we obtain the following.
\begin{thm}\label{KacR formula} For $\lambda\in\mathbb{Z}_+^n$, we have
$${\rm ch}L^0(\Lambda(\lambda))=S^{\A/\B}_{\lambda},$$
where $\A=\mathbb{Z}_{>0}$ and $\B=\mathbb{Z}_{\leq 0}$.\qed
\end{thm}

\begin{cor}
Under the above hypothesis, we have
$${\rm
ch}L^0(\Lambda(\lambda))={\rm det}\left( {\rm
ch}L^0(\Lambda(\lambda_i-i+j)) \right)_{1\leq i,j\leq n}.$$\qed
\end{cor}
\vskip 3mm

\subsection{Grothendieck rings}
Let us summarize the previous results in terms of categories of
representations.

For $n\geq 1$, let $\R_n$ be the Grothendieck group for the category
of rational representations of  $\frak{gl}_n$. We may view $\R_n$
($n\geq 1$) as the free $\mathbb{Z}$-module spanned by the rational
Schur polynomials $s_{\lambda}=s_{\lambda}({\bf x}_{[n]})$ for
$\lambda\in\mathbb{Z}_+^n$. Let $\R=\bigoplus_{n\geq 0}\R_n$, where
$\R_0=\mathbb{Z}$. If we identify the elements in $\R\otimes \R$
with the functions in two sets of variables ${\bf x}_{\mathbb{N}}$
and ${\bf y}_{\mathbb{N}}$, then there exists a natural
comultiplication $\Delta : \R \rightarrow \R\otimes \R $ defined by
$\Delta(1)=1\otimes 1$ and $\Delta(s_{\lambda}({\bf
x}_{[n]}))=\sum_{p+q=n}s_{\lambda}({\bf x}_{[p]},{\bf y}_{[q]})$ for
$\lambda\in\mathbb{Z}_+^n$, where $1$ denotes the unity in $\R_0$
and ${\bf y}_{[q]}=\{\,x_{p+1},\ldots,x_n\,\}$. By \eqref{LR' rule
for rational Schur}, we have
$$\Delta(s_{\lambda}({\bf
x}_{[n]}))=s_{\lambda}({\bf x}_{[n]})\otimes 1 + 1\otimes
s_{\lambda}({\bf x}_{[n]})+
\sum_{\mu,\nu}\hat{c}^{\lambda}_{\mu\,\nu}s_{\mu}({\bf
x}_{[p]})\otimes s_{\nu}({\bf x}_{[q]}),$$ where the sum is taken
over all $\mu\in\mathbb{Z}_+^p$ and $\nu\in\mathbb{Z}_+^q$ such that
$p+q=n$. Also, the counit $\varepsilon : \R \rightarrow \mathbb{Z}$
is the $\mathbb{Z}$-linear map which vanishes on $\R_n$ for $n\geq
1$ with $\varepsilon(1)=1$. Note that $\Delta$ preserves the
grading, that is, $\Delta(\R_n)\subset (\R\otimes
\R)_n=\bigoplus_{p+q=n}\R_p\otimes \R_q$. Hence $\R$ is a graded
cocommutative coalgebra over $\mathbb{Z}$ (cf.\cite{Z}).

Let $\R^*=\bigoplus_{n\geq 0}\R^*_n$ be the graded dual of $\R$,
where $\R^*_n={\rm Hom}_{\mathbb{Z}}(\R_n,\mathbb{Z})$. Let
$\sigma_{\lambda}$ be the element in $\R^*$ dual to $s_{\lambda}$
for $\lambda\in\mathbb{Z}_+^n$ (that is,
$\sigma_{\lambda}(s_{\mu})=\delta_{n\, n'}\delta_{\lambda\,\mu}$ for
$\mu\in\mathbb{Z}_+^{n'}$). Then for any $\varphi \in \R^*$, we may
write
$\varphi=\varphi(1)\varepsilon+\sum_{\lambda}\varphi(s_{\lambda})\sigma_{\lambda}$.

Since $\R$ is a cocommutative coalgebra  over $\mathbb{Z}$, $\R^*$
naturally becomes a graded commutative $\mathbb{Z}$-algebra with the
multiplication $\Delta^*$   given by
\begin{equation}\label{mul of dual R}
\begin{split}
&\Delta^*(\sum_{\mu}a_{\mu}\sigma_{\mu}\otimes
\sum_{\nu}b_{\nu}\sigma_{\nu})=\sum_{\lambda}\left(\sum_{\mu,\nu}a_{\mu}b_{\nu}\hat{c}^{\lambda}_{\mu\,\nu}\right)\sigma_{\lambda},\\
& \Delta^*(\sum_{\mu}a_{\mu}\sigma_{\mu}\otimes \varepsilon)=
\Delta^*(\varepsilon\otimes\sum_{\mu}a_{\mu}\sigma_{\mu})=\sum_{\mu}a_{\mu}\sigma_{\mu},\\
&\Delta^*(\varepsilon\otimes\varepsilon)=\varepsilon,
\end{split}
\end{equation}
where we assume that $\hat{c}^{\lambda}_{\mu\,\nu}=0$ unless
$\lambda\in\mathbb{Z}_+^{p+q}$, $\mu\in\mathbb{Z}_+^p$ and
$\nu\in\mathbb{Z}_+^q$.\vskip 5mm

From now on, we assume that
$\frak{g}=\widehat{\frak{gl}}_{\infty|\infty}$ or
$\widehat{\frak{gl}}_{\infty}$. For each $\lambda\in\mathbb{Z}_+^n$,
we denote by $\mathscr{L}_{\lambda}$ the associated irreducible
representation of $\frak{g}$ (that is,
$\mathscr{L}_{\lambda}=L(\Lambda(\lambda))$ or
$L^0(\Lambda(\lambda))$). For $n\geq 1$, let $\mathcal{O}_n$ be the
category of $\frak{g}$-modules $V$ which are isomorphic to
\begin{equation*}
\bigoplus_{\lambda\in\mathbb{Z}_+^n}\mathscr{L}_{\lambda}^{\oplus
m_{\lambda}}
\end{equation*}
for some $m_{\lambda}\in\mathbb{Z}_{\geq 0}$. For convenience, we
denote by $\mathcal{O}_0$ the category of finite dimensional trivial
representations of $\frak{g}$.

\begin{lem}\label{product of Gr ring}
For $V\in \mathcal{O}_m$ and $W\in\mathcal{O}_{n}$, we have
$V\otimes W\in\mathcal{O}_{m+n}$.
\end{lem}
\pf We may assume that $m,n>0$. Suppose that
\begin{equation*}
V=\bigoplus_{\mu\in\mathbb{Z}_+^m}\mathscr{L}_{\mu}^{\oplus
a_{\mu}}, \ \
W=\bigoplus_{\nu\in\mathbb{Z}_+^n}\mathscr{L}_{\nu}^{\oplus
b_{\nu}},
\end{equation*}
for some $a_{\mu}$ and $b_{\nu}\in\mathbb{Z}_{\geq 0}$. From
\eqref{LR' rule for rational Schur}, we see that for each
$\lambda\in\mathbb{Z}_+^{m+n}$, there exist only finitely many $\mu$
and $\nu$'s such that $\hat{c}^{\lambda}_{\mu\,\nu}\neq 0$. By
Theorem \ref{LR rule for SAB}, the multiplicity of
$\mathscr{L}_{\lambda}$ ($\lambda\in\mathbb{Z}_+^{m+n}$) in
$V\otimes W$ is given by
\begin{equation*}
\sum_{\mu,\nu}a_{\mu}b_{\nu}\hat{c}^{\lambda}_{\mu\,\nu},
\end{equation*}
which is a non-negative integer, and hence $V\otimes
W\in\mathcal{O}_{m+n}$. \qed\vskip 3mm

Let
\begin{equation}
\mathbf{K}(\mathfrak{g})=\bigoplus_{n\geq
0}\mathbf{K}(\mathcal{O}_n)
\end{equation}
be the direct sum of the Grothendieck groups of $\mathcal{O}_n$.
By Lemma \ref{product of Gr ring}, $\mathbf{K}(\mathfrak{g})$
naturally becomes a commutative $\mathbb{Z}$-algebra with the
multiplication given as a tensor product. We denote by $[V]\in
\mathbf{K}(\mathcal{O}_n)$ the isomorphism class of
$V\in\mathcal{O}_n$.

For $n\geq 1$, we define a $\mathbb{Z}$-linear map $\chi_n : \R^*_n
\rightarrow \mathbf{K}(\mathcal{O}_n)$ by
\begin{equation}
\chi_n(\varphi)=\sum_{\lambda\in\mathbb{Z}_+^n}\varphi(s_{\lambda})\left[\mathscr{L}_{\lambda}
\right],
\end{equation}
for $\varphi\in\R^*_n$. Also, we define
$\chi_0(\varepsilon)=\left[\mathbb{C}\right]$, the isomorphism class
of the one-dimensional trivial representation. So we have a
$\mathbb{Z}$-linear map
\begin{equation}
\chi=\bigoplus_{n\geq 0}\chi_n : \R^* \longrightarrow
\mathbf{K}(\mathfrak{g}).
\end{equation}

\begin{thm}\label{character ring}\mbox{}
$\chi$ is an isomorphism of commutative $\mathbb{Z}$-algebras.
\end{thm}
\pf It follows from Theorem \ref{CL formula}, \ref{KacR formula},
and \ref{LR rule for SAB}. \qed\vskip 3mm

Let us end this section with a remark on an involution on $\R^*$. We
define a $\mathbb{Z}$-linear map $\omega : \R^* \rightarrow \R^*$ by
\begin{equation}
\omega(\sum_{\lambda}a_{\lambda}\sigma_{\lambda})=\sum_{\lambda}a_{\lambda}\sigma_{\lambda^*},
\end{equation}
for $\sum_{\lambda}a_{\lambda}\sigma_{\lambda}\in\R^*$, and
$\omega(\varepsilon)=\varepsilon$. It is clear that $\omega$ is an
isomorphism of $\mathbb{Z}$-modules.

\begin{lem}\label{symmetry of LR}\mbox{}
There exists a bijection from ${\bf LR}^{\lambda}_{\mu\,\nu}$ to
${\bf LR}^{\lambda^*}_{\mu^*\,\nu^*}$ for
$\lambda\in\mathbb{Z}_+^{m+n},\mu\in\mathbb{Z}_+^m$, and
$\nu\in\mathbb{Z}_+^n$.
\end{lem}
\pf Given $\lambda\in\mathbb{Z}_+^n$, let $p,q$ be sufficiently
large positive integers such that
$\lambda+(p^n),\lambda^*+(q^n)\in\cP_n$. Then we can check that
\begin{equation*}
\delta^{p+q}_n([\lambda+(p^n)]')=[\lambda^*+(q^n)]',
\end{equation*}
as Young diagrams.

Now, suppose that
$\lambda\in\mathbb{Z}_+^{m+n},\mu\in\mathbb{Z}_+^m,
\nu\in\mathbb{Z}_+^n$ are given, and $p,q$ are sufficiently large
positive integers. Then, the required bijection comes from the
following one-to-one correspondences
\begin{equation*}
\begin{split}
{\bf LR}^{\lambda}_{\mu\,\nu}&\stackrel{1-1}{\longleftrightarrow} LR^{[\lambda+(p^{m+n})]'}_{[\mu+(p^m)]'\,[\nu+(p^n)]'}\\
&\stackrel{1-1}{\longleftrightarrow}LR^{\delta^{p+q}_{m+n}([\lambda+(p^{m+n})]')}_{\delta^{p+q}_{m}([\mu+(p^m)]')\,\delta^{p+q}_{n}([\nu+(p^n)]')}
\ \ \text{by Lemma \ref{dual LR}} \\ & \ \ \ \
=LR^{[\lambda^*+(q^{m+n})]'}_{[\mu^*+(q^m)]'\,[\nu^*+(q^n)]'}
\\
&\stackrel{1-1}{\longleftrightarrow}{\bf
LR}^{\lambda^*}_{\mu^*\,\nu^*}.
\end{split}
\end{equation*}
\qed

\begin{rem}{\rm
The map $\theta : \R \rightarrow \R$ sending $s_{\lambda}({\bf
x}_{[n]})$ to $s_{\lambda}({\bf x}_{[n]}^{-1})$ for
$\lambda\in\mathbb{Z}_+^n$ with $\theta(1)=1$, is an automorphism of
a coalgebra over $\mathbb{Z}$. In fact, $\theta(s_{\lambda}({\bf
x}_{[n]}))=s_{\lambda^*}({\bf x}_{[n]})$, and $\omega$ is the map
induced from $\theta$. From \eqref{LR' rule for rational Schur} and
the linear independence of rational Schur polynomials, it follows
directly that
$\hat{c}^{\lambda}_{\mu\,\nu}=\hat{c}^{\lambda^*}_{\mu^*\,\nu^*}$
for $\lambda\in\mathbb{Z}_+^{m+n},\mu\in\mathbb{Z}_+^m,
\nu\in\mathbb{Z}_+^n$, while Lemma \ref{symmetry of LR} gives a
bijective proof  of this. }
\end{rem}

Therefore, it follows that
\begin{prop}
$\omega$ is an involution of $\R^*$ as a $\mathbb{Z}$-algebra.\qed
\end{prop}
\end{document}